\documentclass[pdflatex,sn-mathphys-num]{sn-jnl}

\usepackage{amsmath,amssymb,amsfonts}
\let\orcidlogo\undefined
\usepackage{orcidlink}
\usepackage{amsthm}
\usepackage{microtype}
\usepackage{mathtools}
\usepackage{comment}
\usepackage{enumitem}
\usepackage[title]{appendix}
\hypersetup{colorlinks=true,citecolor=blue,linkcolor=blue,urlcolor=blue}
\numberwithin{equation}{section}

\theoremstyle{thmstyleone}
\newtheorem{theorem}{Theorem}[section]
\newtheorem{proposition}[theorem]{Proposition}
\newtheorem{lemma}[theorem]{Lemma}
\newtheorem{corollary}[theorem]{Corollary}
\newtheorem{definition}[theorem]{Definition}
\newtheorem{convention}[theorem]{Notation}
\theoremstyle{thmstyletwo}
\newtheorem{remark}[theorem]{Remark}
\newcommand{\leanrepo}{https://github.com/itpplasma/stafford38-formal/blob/1029fdb382eee06e7b59bcf32801f40931e2079a}
\newcommand{\leanlibrepo}{https://github.com/itpplasma/algebraic-analysis/blob/bbbbf3fc358ca8100b158cec4cf47f336ab70163}
\newcommand{\leanglobalrepo}{https://github.com/itpplasma/global-stafford-formal/blob/e676caed1253e5d09fc6699d5dc683febbabf7bb}
\newcommand{\leanfilteredrepo}{https://github.com/itpplasma/symplectic-stafford-formal/blob/dccff0d8fe723ec1478d1f63bbbdf648ba240bfa}
\newcommand{\leandecl}[3]{\href{\leanrepo/#1\#L#2}{\nolinkurl{#3}}}
\newcommand{\leanlib}[3]{\href{\leanlibrepo/#1\#L#2}{\nolinkurl{#3}}}
\newcommand{\leanglobal}[3]{\href{\leanglobalrepo/#1\#L#2}{\nolinkurl{#3}}}
\newcommand{\leanfiltered}[3]{\href{\leanfilteredrepo/#1\#L#2}{\nolinkurl{#3}}}
\newif\ifshowleandetails
\showleandetailsfalse
\newcommand{\formalref}[2][]{
  \ifshowleandetails
    {\small\textsf{Lean}\if\relax\detokenize{#1}\relax\else\ (#1)\fi: #2}%
  \fi}

\newcommand{\Char}{\operatorname{Char}}

\newcommand{\ord}{\operatorname{ord}}
\newcommand{\gr}{\operatorname{gr}}
\newcommand{\Ann}{\operatorname{Ann}}

\newcommand{\Spec}{\operatorname{Spec}}
\newcommand{\Supp}{\operatorname{Supp}}

\begin{document}
 \setlength{\parindent}{0pt}
\title[Stafford's conjecture on cyclicity of torsion modules]{Stafford's conjecture on cyclicity of torsion modules}

\author[1]{\fnm{Christopher} \sur{Albert}\,\orcidlink{0000-0003-4773-416X}}\email{albert@tugraz.at}

\author*[1]{\fnm{Johanna} \sur{Moser}\,\orcidlink{0009-0005-5289-6915}}\email{j.moser@tugraz.at}

\author[1]{\fnm{Maximilian} \sur{Philipp}\,\orcidlink{0009-0009-0379-7637}}\email{philipp@student.tugraz.at}

\affil[1]{\orgname{Graz University of Technology}, \orgaddress{\city{Graz},
  \country{Austria}}}

\abstract{%
We prove that over a Weyl algebra over a field of characteristic zero $A_n(k)$, every nonzero $d\in A_n(k)$ admits $F,R,S\in A_n(k)$ such that
$1=dR+FdS$, a conjecture Stafford proposed in 1978. As a consequence, every finitely generated torsion right $A_n(k)$-module is cyclic.
The proof maps $d$, by a symplectic change of variables, to a normal form $\tilde d$
and shows that the cyclic module $Q=A_n(k)/(\tilde dA_n(k)+x_1^N\tilde dA_n(k))$ vanishes.
We first show that right multiplication by $x_1$ is surjective on $Q$. A comparison of lengths of localized modules
over the associated graded ring then gives $\operatorname{gr}Q=x_1\operatorname{gr}Q$,
so that the characteristic variety of $Q$ is disjoint from the hyperplane
$\{x_1=0\}$. Finally, Gabber's involutivity theorem together with an argument about conormal directions at infinity from algebraic geometry shows that such a characteristic variety must be empty and therefore, $Q = 0$.
}



\maketitle

\section{Introduction}

Let $k$ be a field of characteristic zero, and let
\[ A_n(k)=k\langle x_1,\ldots,x_n,\partial_1,\ldots,\partial_n\rangle\]
be the $n$th Weyl algebra, with $[\partial_i,x_j]=\delta_{ij}$ and
$[\partial_i,\partial_j]=[x_i,x_j]=0$.

Stafford's Conjecture~3.8 in \cite{Stafford1978} asks whether every nonzero $d\in A_n(k)$ admits $F\in A_n(k)$ with $A_n(k)=dA_n(k)+FdA_n(k)$
\cite[p.~438]{Stafford1978}. Stafford shows that this formulation implies stable cyclicity for all finitely generated torsion right Weyl-algebra modules
\cite[pp.~437--438]{Stafford1978}. 
We prove Stafford's conjecture:
\begin{theorem}\label{thm:main}
For every field of characteristic zero $k$, every $n\geq0$, and every
nonzero $d\in A_n(k)$, there exists an $F\in A_n(k)$ such that
\[ A_n(k)=dA_n(k)+FdA_n(k),\]
or equivalently, there exist $F, S, R\in A_n(k)$ such that
\[ 1=dR+FdS.\]
Moreover, if $n\ge1$, one may choose $F=\Phi_g^{-1}(x_1)^{\deg_Bd}$ for a linear
$\omega$-preserving automorphism $g$ as in Section~\ref{sec:monic}. 
\end{theorem}

\formalref{\leandecl{Stafford38/Paper/Main.lean}{24}{Stafford38.Paper.theorem_1_1} proves \leandecl{Stafford38/ChallengeDefinitions.lean}{78}{Stafford38Challenge.UniversalStatement};
\leandecl{Stafford38/Paper/Main.lean}{29}{Stafford38.Paper.theorem_1_1_fixedSource} proves \leandecl{Stafford38/ChallengeDefinitions.lean}{193}{Stafford38FixedSourceChallenge.UniversalFixedSourceStatement}
with \leandecl{Stafford38/ChallengeDefinitions.lean}{170}{Stafford38FixedSourceChallenge.IsLinearWeylCoordinate} and \leandecl{Stafford38/ChallengeDefinitions.lean}{155}{Stafford38FixedSourceChallenge.bernsteinDegree};
Mathlib-only statements compared by Comparator: \leandecl{Stafford38/ChallengeDefinitions.lean}{291}{Stafford38Challenge.PaperTheorem_1_1},
\leandecl{Stafford38/ChallengeDefinitions.lean}{78}{Stafford38Challenge.UniversalStatement},
\leandecl{Stafford38/ChallengeDefinitions.lean}{193}{Stafford38FixedSourceChallenge.UniversalFixedSourceStatement}.}

\textit{Proof sketch:}
We state some preliminaries in section \ref{sec:setting}, then start in section \ref{sec:monic} by defining a map $\Phi_g$ sending any chosen $d$ of Bernstein degree $N$ to a monic normal form $\tilde d$ with order symbol $P$ using Bernstein and order filtrations of $A_n(k)$. For $\tilde d$, we prove in sections \ref{sec:euler}-\ref{sec:proof} that the module $Q = A_n(k)/(\tilde d A_n(k) + x_1^N\tilde d A_n(k)) $ vanishes, and find the final form of Theorem \ref{thm:main} by applying the inverse map $\Phi_g^{-1}$ to $\tilde d A_n(k) + x_1^N\tilde d A_n(k).$ The main components of the proof are the following: in section \ref{sec:euler}, we show that $Q = Qx_1$ using the extended euclidean algorithm, and $\gr Q = x_1\gr Q,$ through two contradictory length estimations of localized modules. In section \ref{sec:conormal-review}, using mainly results from algebraic geometry, we show that a set $W$ fulfilling certain criteria (Zariski closed over an algebraically closed field of characteristic zero, fibre-conical, involutive, there exists a polynomial with $P(1,0\dots,0)\neq 0$ vanishing on $W$, and $W$ does not contain points with $x_1 = 0$) must always be empty (Theorem \ref{thm:coisotropic-exclusion}). Finally, we show in section \ref{sec:proof} that by defining $W$ to be the characteristic variety of the closed-field pendant to $Q$, the results of section \ref{sec:euler} as well as the definition of $P$ in section \ref{sec:monic} ensure that $W$ fulfills all criteria and must thus be empty. We conclude that therefore, also $Q$ itself must be zero, from which Theorem \ref{thm:main} follows directly. Lastly, we provide the proof by Stafford \cite{Stafford1978} that hence, every finitely generated torsion module is stably cyclic (Corollary \ref{cor:cyclic}).

\subsection{Related literature}
Dixmier showed that over the first Weyl algebra $A_1(k)$ with $k$ having characteristic zero, any right ideal can be generated by at most two elements \cite{dixmier}. In a seminal paper in 1978, Stafford extended this result to all $A_n(k)$ and moreover showed, that any finitely generated left torsion $A_n(k)$-module is a homomorphic image of a projective left ideal, and as such can also be (stably) generated by at most two elements \cite[Theorem 3.7]{Stafford1978}.
Hillebrand and Schmale \cite{HillebrandSchmale2001}, as well as Leykin \cite{Leykin2004} give an algorithm to find the two generators of left ideals, also valid for rational coefficient algebras $B_n(k)$. Caro and Levcovitz \cite{CaroLevcovitz2010} generalize the result for both left and right ideals to formal and to convergent Laurent series, and Quadrat and Robertz \cite{QuadratRobertz2014} provide a constructive proof method to calculate the two generators of torsion modules over the Weyl algebras, as well as rings with rational coefficients, formal or convergent Laurent series, and, for ordinary differential operators, with formal or locally convergent power series coefficients. 

In the same paper \cite{Stafford1978}, Stafford states the Conjecture 3.8 (our Theorem \ref{thm:main}), and shows that this conjecture implies that any finitely generated right torsion module can be stably generated by one element, i.e. is stably cyclic (Corollary \ref{cor:cyclic}). He shows that cyclicity is true if $n = 1$ (since $A_1(k)$ has the restricted minimum condition which makes torsion modules Artinian and he proves that Artinian modules are cyclic) and also states that the conjecture is true if $d$ is monomial, but provides, as far as the authors are aware, no proof of this. This case is also contained in our proof.
For \textit{holonomic} modules over ``simple rings such that they do not have finite length over themselves as left modules'' (which all Weyl algebras are), Leykin \cite{Leykin2004} proved cyclicity constructively, and attributed the proof method to Stafford. However, torsion modules need not be holonomic: for $n\ge2$, Stafford constructed quotients $A_n(k)/dA_n(k)$ of Gelfand--Kirillov dimension $2n-1$ as counterexamples  \cite{StaffordNonholonomic1985}, so the question about cyclicity remained open for non-holonomic torsion modules. A paper by Bellamy in 2026 states this explicitly (\cite{Bellamy2026} Conjecture 3.4). 

We base our proof on the Weyl algebra and on filtrations and characteristic
varieties of $A_n(k)$-modules as treated in the textbooks
\cite{coutinho,Bjork,HTT,MR,SchnellNotes}, on the involutivity of
the characteristic variety due to Gabber \cite{Gabber1981,Gabber1982,SinghKumar2014}, on commutative algebra, in particular localization, length, valuation rings, Kähler differentials and completions
\cite{atiyah,matsumura,Stacks}, field theory and bilinear forms \cite{lang} and spectral sequences \cite{weibel}, and on
algebraic geometry, mostly based on \cite{hartshorne}.

\subsection{Formal verification with Lean}\label{sec:formal}
Theorem~\ref{thm:main} and
Corollary~\ref{cor:cyclic} are proved in Lean~4~\cite{Lean4} with
Mathlib~\cite{Mathlib2020} and AlgebraicAnalysis
v0.3.3~\cite{AlgebraicAnalysis033}. The proofs use only Lean's three
standard axioms: propositional extensionality (\texttt{propext}), the axiom
of choice (\texttt{Classical.choice}) and soundness of quotient types
(\texttt{Quot.sound}). Every intermediate result is proved in Lean as well,
including Gabber's involutivity theorem for cyclic $A_n(k)$-modules. The
proof in Sections~\ref{sec:monic}--\ref{sec:proof} follows the formal proof
step by step. The main theorem is registered
in Palomar~\cite{Stafford38Palomar2026}, and the Lean sources are
archived~\cite{Stafford38Zenodo140}. Tooling and the role of AI assistance
are described in Appendix~\ref{app:ai-workflow}. Section~\ref{sec:generalizations} states further consequences proved in Lean and two generalizations, to rings of differential operators on smooth affine varieties and to simple filtered algebras.
The supplementary
material~\cite{Stafford38Supplementary033} shows each formalized statement of
this paper next to its Lean version and proof. It can be browsed at
\url{https://itpplasma.github.io/stafford38-supplementary/}.

\section{Preliminaries}\label{sec:setting}

Throughout, $k$ is a field of characteristic zero and for the $n$-th Weyl algebra over $k$, we write
\[ A_n(k)=k\langle x_1,\dots,x_n,\partial_1,\dots,\partial_n\rangle,
  \quad [\partial_i,x_j]=\delta_{ij},\quad [x_i,x_j]=[\partial_i,\partial_j]=0 .\]
All ideals are right ideals and all modules are right modules
unless the contrary is said explicitly.
\begin{definition}[e.g. \cite{coutinho}]\label{def:cyclic}
Let $M$ be a right $A_n(k)$-module. We say $M$ is \emph{torsion} if for every $m\in M$ there exists
$0\ne d\in A_n(k)$ with $md=0$, and \emph{cyclic} if $M=mA_n(k)$ for some $m$. 
\end{definition}
\begin{convention}\label{conv:pair}
We single out the first Weyl pair $(x_1, \partial_1)$ and write
\[\qquad p:=\partial_1,\qquad
  C:=A_{n-1}(k)[x_1],\qquad A_n(k)=C[p;\delta],\]
Here $\delta(x_1)=1$ and $\delta$ vanishes on $A_{n-1}(k)$,
so that $x_1$ is central in $C$, $[p,x_1] = 1$, and the \emph{transverse algebra} $A_{n-1}(k)=k\langle x_i,\partial_i : i\ge2\rangle$ commutes with $x_1$ and $p$.
Every element of $A_n(k)$ can uniquely be written as $\sum_j p^j c_j$ with $c_j\in C$. \cite[Prop. 1.2.1]{coutinho}
\end{convention}
\formalref{%
\leandecl{Stafford38/Paper/Section2/Weyl.lean}{98}{Stafford38.Paper.Section2.transverse},
\leandecl{Stafford38/Paper/Section2/Weyl.lean}{103}{Stafford38.Paper.Section2.C},
\leandecl{Stafford38/Paper/Section2/Weyl.lean}{179}{Stafford38.Paper.Section2.convention_2_2_relations},
\leandecl{Stafford38/Paper/Section2/Weyl.lean}{199}{Stafford38.Paper.Section2.convention_2_2_exists},
\leandecl{Stafford38/Paper/Section2/Weyl.lean}{228}{Stafford38.Paper.Section2.convention_2_2_unique}}

\begin{definition}[{[Def.~2.6 and ~2.9 \cite{SchnellNotes}]}]\label{def:filtration}
A filtration $F_\bullet$ on $A_n(k)$ is a sequence of linear subspaces
\[\{0\} = F_{-1} \subseteq F_0 \subseteq F_1\subseteq...\subseteq A_n(k),\]
such that $F_i\cdot F_j \subseteq F_{i+j}$ and $A_n(k) = \bigcup F_j.$
The associated graded algebra, is defined as
\[\mathrm{gr}^F A_n(k) = \bigoplus_{j = 0}^\infty F_j/F_{j-1}.\]
Let $a\neq 0 \in A_n(k)$ and let $j$ be the smallest integer, such that $a \in F_j$, called a \emph{filtration degree} $\deg_F a$. The class of $a$ in $F_j/F_{j-1}$ is called its \emph{symbol} $\sigma_j(a)$. A symbol, being the class of an element in a single graded piece, is automatically homogeneous of the corresponding
degree, i.e.\ contains only terms of that total degree.
\end{definition}

In the proof, we use two filtrations of $A_n(k)$ (for more details, see e.g. \cite[Chapter 7]{coutinho})
\begin{itemize}[nosep,leftmargin=1.4em]
\item The \emph{Bernstein filtration} $B_\bullet$ assigns degree $1$ to every
  $x_i$ and every $\partial_i$. The Bernstein degree of $d$ is
  written $\deg_B d$.
\item The \emph{order filtration} $F_\bullet$ assigns degree $0$ to every $x_i$
  and degree $1$ to every $\partial_i$; $F_m$ is spanned by the monomials
  $x^\alpha\partial^\beta$ with $|\beta|\le m$.  
\end{itemize}
\formalref{%
\leandecl{Stafford38/Paper/Section2/Filtration.lean}{79}{Stafford38.Paper.Section2.IsFiltration},
\leandecl{Stafford38/Paper/Section2/Filtration.lean}{68}{Stafford38.Paper.Section2.gr},
\leandecl{Stafford38/Paper/Section2/Filtration.lean}{86}{Stafford38.Paper.Section2.filtrationDegree},
\leandecl{Stafford38/Paper/Section2/Filtration.lean}{95}{Stafford38.Paper.Section2.symbol},
\leandecl{Stafford38/Paper/Section2/Filtration.lean}{104}{Stafford38.Paper.Section2.bernsteinFiltration},
\leandecl{Stafford38/Paper/Section2/Filtration.lean}{109}{Stafford38.Paper.Section2.orderFiltration},
\leandecl{Stafford38/Paper/Section2/Filtration.lean}{121}{Stafford38.Paper.Section2.bernsteinDegree},
\leandecl{Stafford38/Paper/Section2/Filtration.lean}{185}{Stafford38.Paper.Section2.grEquivSymbolRing}}
For both the order filtration and the Bernstein filtration, the associated graded algebra $\mathrm{gr} A_n(k)$ is simply the polynomial ring $k[x,\xi]=k[x_1,\dots,x_n,\xi_1,\dots,\xi_n]$ in $2n$ variables with $\xi_i$ being the symbol of $\partial_i$ \cite[Th. 7.3.1]{coutinho}. In particular, $\mathrm{gr} A_n(k)$ is commutative \cite[Prop.~2.10]{SchnellNotes}.
Good filtrations are treated systematically for example in \cite{Bjork,HTT}. 

\begin{definition}[{\cite[Ch. 8.3]{coutinho}}]\label{def:good-filtration}
A filtration of a right $A_n$-module $M$ compatible with a filtration $F_\bullet$
is an exhaustive increasing chain of $k$-subspaces $M_\bullet$ with
$M_iF_j\subseteq M_{i+j}$ and $M_i=0$ for $i<0$. The associated graded module is $\operatorname{gr}M=\bigoplus_iM_i/M_{i-1}$; it is a graded module over $\operatorname{gr}_FA_n(k)$ via the maps
$M_i/M_{i-1}\times F_j/F_{j-1}\to M_{i+j}/M_{i+j-1}$, $(\bar m,\bar a)\mapsto\overline{ma}$,
which are well defined because $M_{i-1}F_j\subseteq M_{i+j-1}$ and $M_iF_{j-1}\subseteq M_{i+j-1}$.
It is \emph{good} if the associated graded module is finitely generated over $\gr^FA_n$.
\end{definition}
As an example, the filtration $M_i:=qF_i$ of a cyclic module over $A_n$ induced by its generator $q$ is good, since $\gr Q$ is generated by the element $\bar{q}$ denoting the class of $q$ in $Q_0/Q_{-1}.$

For a finitely generated right $A_n$-module $M$ with a good order filtration
$M_\bullet$, the graded module $\gr M$ is a finitely generated $k[x,\xi]$-module.

\begin{definition}[{\cite[Ch. 1]{atiyah}}]\label{def:spec}
Let $S$ be a commutative ring. The \emph{spectrum} $\operatorname{Spec}S$ is the set
of all prime ideals of $S$. We say that $f\in S$ vanishes at the point
$\mathfrak p\in\operatorname{Spec}S$ if $f\in\mathfrak p$. The \emph{vanishing locus}
of an ideal $J\subseteq S$ is
\[ V(J)=\{\mathfrak p\in\operatorname{Spec}S:\ \mathfrak p\supseteq J\},\]
the set of points at which all elements of $J$ vanish; these sets are the (Zariski) closed sets
of the Zariski topology on $\operatorname{Spec}S$. For a field $K$ and an ideal $J\subseteq K[z_1,\dots,z_N]$ we write
\[V_K(J)=\{a\in K^N : f(a)=0\ \forall f\in J\}\]
for its (classical) zero set, and for $S\subseteq K^N$ we write
$I(S)=\{f\in K[z]: f(a)=0\ \forall a\in S\}$ for the \emph{vanishing ideal} of $S$.
The sets $V_K(J)$ are the closed sets of the Zariski topology on $K^N$.
Since for $f$ in the maximal ideal $f\in\mathfrak m_a:=(z_1-a_1,\dots,z_N-a_N)$ iff $f(a)=0$, we have
$a\in V_K(J)$ if and only if $\mathfrak m_a\in V(J)$.\\
For an ideal $J$ of $S$, its \emph{radical} is
$\sqrt J=\{f\in S:\ f^m\in J\text{ for some }m\ge1\}$. $J$ is called \emph{radical} if $J=\sqrt J$. Since a prime ideal contains $J$
iff it contains $\sqrt J$, we have $V(J)=V(\sqrt J)$.
\end{definition}

\begin{definition}[{\cite[Appendix E]{HTT}}]\label{def:cotangent}
Let $\mathbb A^n_k$ be affine $n$-space over $k$ with coordinates $x_1,\dots,x_n$.
Its \emph{cotangent bundle} $T^*\mathbb A^n_k=\mathbb A^{2n}_k$ has the additional
fibre coordinates $\xi_1,\dots,\xi_n$ dual to $x_1,\dots,x_n$, so that a covector is
$\sum_i\xi_i\,dx_i$. Its coordinate ring is $k[x,\xi]=\gr A_n(k)$, and we regard
$\mathbb A^n_k=\operatorname{Spec}k[x]$ and $T^*\mathbb A^n_k=\operatorname{Spec}k[x,\xi]$
as sets of prime ideals \cite[Tag 00DZ]{Stacks}.
\end{definition}

\formalref{\leandecl{Stafford38/Paper/Section2/Geometry.lean}{159}{Stafford38.Paper.Section2.CotangentBundle},
\leandecl{Stafford38/Characteristic/Polynomial.lean}{22}{Stafford38.Characteristic.SymbolRing}
(the ring $k[x,\xi]$; the first summand indexes base coordinates $x_i$, and the second indexes fibre coordinates $\xi_i$)}

\begin{definition}[{\cite[Ch. 2.2]{HTT}}]\label{def:char-variety}
    The \emph{characteristic variety} of a finitely generated $A_n$-module $M$ with a good order filtration
$M_\bullet$ is defined by
\[\Char(M):=V\bigl(\sqrt{\Ann_{k[x,\xi]}\gr M}\bigr) = V\bigl(\Ann_{k[x,\xi]}\gr M\bigr)\subseteq
  T^*\mathbb A^n_k=\Spec k[x,\xi],\]
$\Ann_{k[x,\xi]}\gr M = \{f \in k[x,\xi]: f\cdot m = 0\,\forall m \in \gr M\}$ are the \emph{annihilators} of $\gr M$ over $k[x,\xi]$, and $V$ denotes their vanishing locus.

$\Char(M)$ is independent of the choice of the good filtration. \cite[Lemma D.3.1]{HTT}
\end{definition}
\formalref{%
\leandecl{Stafford38/Paper/Section2/CharacteristicVariety.lean}{54}{Stafford38.Paper.Section2.CharOf},
\leandecl{Stafford38/Paper/Section2/CharacteristicVariety.lean}{83}{Stafford38.Paper.Section2.Char},
\leandecl{Stafford38/Paper/Section2/CharacteristicVariety.lean}{61}{Stafford38.Paper.Section2.definition_2_7_independent},
\leandecl{Stafford38/Paper/Section2/CharacteristicVariety.lean}{87}{Stafford38.Paper.Section2.Char_eq}}

\begin{definition}[{\cite[Ch. E.1]{HTT}}]\label{def:symplectic}
A bilinear form $\omega$ on a $k$-vector space is \emph{symplectic} if it is
antisymmetric and nondegenerate. The \emph{standard symplectic form} on $k^{2n}$
with basis $x_1,\dots,x_n,\partial_1,\dots,\partial_n$ is determined by
\[ \omega(x_i,\partial_j)=\delta_{ij},\qquad
 \omega(x_i,x_j)=\omega(\partial_i,\partial_j)=0,\]
analogously to the Weyl relations. 
The form is \emph{nondegenerate} if $\omega(u,\cdot)=0$ implies $u=0$. For a subspace $U$, we define 
$U^\perp=\{v:\omega(u,v)=0\ \forall u\in U\}$. A \emph{symplectic basis} is a basis
$e_1,\dots,e_n,f_1,\dots,f_n$ with $\omega(e_i,f_j)=\delta_{ij}$ and
$\omega(e_i,e_j)=\omega(f_i,f_j)=0$. The \emph{Gram matrix} of $\omega$ in a basis
$(b_i)$ is $J=(\omega(b_i,b_j))_{ij}$.
\end{definition}

\begin{definition}[{\cite[Appendix D.3]{HTT}}]\label{def:poisson}
The \emph{Poisson bracket} on $k[x,\xi]$ is
\[ \{f,g\}=\sum_{i=1}^n\left(
   \frac{\partial f}{\partial x_i}\frac{\partial g}{\partial \xi_i}
   -\frac{\partial f}{\partial \xi_i}\frac{\partial g}{\partial x_i}\right).\]
An ideal $J$ is \emph{involutive} if $\{J,J\}\subseteq J$. 

\end{definition}

\formalref{%
\leandecl{Stafford38/Paper/Section2/Symplectic.lean}{130}{Stafford38.Paper.Section2.poisson},
\leandecl{Stafford38/Paper/Section2/Symplectic.lean}{137}{Stafford38.Paper.Section2.poisson_eq}}

\formalref{\leandecl{Stafford38/Paper/Section2/Symplectic.lean}{141}{Stafford38.Paper.Section2.IsInvolutive},
\leandecl{Stafford38/Paper/Section2/Symplectic.lean}{152}{Stafford38.Paper.Section2.IsInvolutiveSet},
\leandecl{Stafford38/Paper/Section2/Symplectic.lean}{146}{Stafford38.Paper.Section2.isInvolutive_iff};
negative control \leandecl{Stafford38/Geometry/GeneralCoisotropicSetsTest.lean}{164}{Stafford38.Geometry.GeneralCoisotropicSetsTest.zeroSectionIdealOne_not_isPoisson}}

\begin{theorem}[Gabber {\cite{Gabber1981}\cite{SinghKumar2014}}{\cite[{Thm. D.3.4}]{HTT}}]\label{ithm:gabber}
Let $K$ be a field of characteristic zero and $M$ a finitely generated $A_n(K)$-module
with a good order filtration. Let $J = \sqrt{Ann_{K[x, \xi]}\gr M}$ be the radical of the annihilators of $M$. Then $J$ is involutive, i.e. \[\{J, J\} \subset J.\]

\end{theorem}
\formalref{%
\leandecl{Stafford38/Paper/Section2/Symplectic.lean}{161}{Stafford38.Paper.Section2.theorem_2_11},
\leandecl{Stafford38/Paper/Section2/Symplectic.lean}{174}{Stafford38.Paper.Section2.theorem_2_11_cyclic}}

\begin{theorem}[{Hilbert's Nullstellensatz, see e.g. \cite[Thm. I.1.3A, Cor. I.1.4, Thm. I.3.2]{hartshorne} and \cite[Cor 7.10, Ex. 7.14]{atiyah}}]\label{thm:nss}
Let $K$ be an algebraically closed field and $J\subseteq K[z_1,\dots,z_N]$ an
ideal. Then
\begin{enumerate}
\item[(a)] the maximal ideals of $K[z]$ are exactly $\mathfrak m_a:=(z_1-a_1,\dots,z_N-a_N)$,
  $a\in K^N$; hence $a\mapsto\mathfrak m_a$ identifies $K^N$ with the closed
  points of $\operatorname{Spec}K[z]$, and $V_K(J)$ with the closed points of
  $V(J)$;
\item[(b)] $V_K(J)=\emptyset$ if and only if $1\in J$;
\item[(c)] $I(V_K(J))=\sqrt J$.
\end{enumerate}
\end{theorem}

\section{Mapping to a normal form}\label{sec:monic}
Given a fixed nonzero $d \in A_n(k)$ of Bernstein degree $N$, our first goal is to find a map to replace $d$ by a normal form $\tilde{d}$ that is monic in $p$, i.e. of the form $\tilde d=p^N+\sum_{j<N}p^ja_j$ with coefficients $a_j$ in $C$. \\

If $N=0$, then $d \in B_0 = k$ is a nonzero scalar, hence a unit, and $1 = dR+FdS$ holds trivially with $F=1$, $R=d^{-1}$, $S=0$. We therefore assume $N\geq1$.

We want to construct a well-behaved automorphism of $\gr^BA_n(k)=k[x,\xi]$
from simpler linear algebra. We proceed in three stages:
define an automorphism of the vector space $\mathbb V=k^{2n}$, extend it to
an automorphism of $A_n(k)$ itself, and observe the induced action this
extension leaves on the associated graded ring.

We can identify $\mathbb V = k^{2n}$ with the Bernstein-degree-one part $B_1/B_0$ of
$A_n(k)$, i.e.\ linear combinations of $x_1,\ldots,x_n,\partial_1,\ldots,
\partial_n$; the coordinates of a vector are then the coefficients of
these generators. By this identification, $\mathbb V$ carries the standard symplectic form $\omega$ (Def. \ref{def:symplectic}). 

\begin{lemma}[{cf. \cite[Thm. 1.1]{dasilva}}]\label{lem:sympcompl}
Let $(\mathbb V,\omega)$ be a symplectic $k$-space of dimension $2n$ and $0\neq v\in \mathbb V$.
Then $v$ is the first vector of a symplectic basis of $\mathbb V$. Consequently, there
exists a $k$-linear automorphism $g$ preserving $\omega$ (i.e. $\omega(gu,gv)=\omega(u,v)$ for all
$u,v\in\mathbb V$) with $g(e)=v$ for any prescribed nonzero vector $e$.
\end{lemma}

\begin{proof}
Since $\omega$ is nondegenerate and $v\neq0$, the functional $\omega(v,-)$
is not identically zero, so there exists a $w'\in\mathbb V$ with
$\omega(v,w')=c\neq0$ and defining $w:=w'/c$ gives $\omega(v,w)=1$.

Let $U=\operatorname{span}(v,w)$. For $0\neq av+bw\in U$: if $a\neq0$,
$\omega(av+bw,w)=a\neq0$; if $b\neq0$, $\omega(av+bw,v)=-b\neq0$. Either
way $av+bw$ pairs nontrivially with something in $U$, so $\omega|_U$ is
nondegenerate.

Since $\omega|_U$ is nondegenerate, $U\cap U^\perp=0$, and as
$\dim U+\dim U^\perp=2n = \dim\mathbb V$ \cite[Prop. 11.2.1]{coutinho} for a nondegenerate form, this means
$\mathbb V=U\oplus U^\perp$. If $u\in U^\perp$ pairs to zero with all of
$U^\perp$, it also pairs to zero with $U$ (as $u\in U^\perp$), hence with
all of $\mathbb V$, forcing $u=0$; so $\omega|_{U^\perp}$ is nondegenerate
too.

By induction on $n$ (trivial for $n=1$, where $U=\mathbb V$), a
symplectic basis $e_2,\ldots,f_n$ of $U^\perp$ therefore exists \cite[XV.8.1 and XV.8.2]{lang}. Setting $e_1:=v$,
$f_1:=w$ gives a symplectic basis of $\mathbb V=U\oplus U^\perp$, since
cross-pairings between $U$ and $U^\perp$ vanish by definition of
$U^\perp$.

Consequently: given nonzero $v,e\in\mathbb V$, we can complete both vectors to a symplectic basis as above. We then define the linear map $g$ to send $e$ to $v$, and each subsequent $e$-basis vector to the corresponding $v$-basis vector, thus preserving all pairwise values of $\omega$ on basis vectors. Through the bilinearity of $\omega$, $g$ then preserves $\omega$ on all of $\mathbb V$; and $g(e)=v$.
\end{proof}
\formalref{%
\leandecl{Stafford38/Paper/Section3/Symplectic.lean}{269}{Stafford38.Paper.Section3.lemma_3_1}}
Let $g$ be a $k$-linear automorphism of $\mathbb V$ that preserves $\omega$.
Since $\mathbb V$ is identified with the span of the generators of $A_n(k)$,
$g$ assigns to each generator $u$ a linear combination $g(u)$ of generators.
Because $A_n(k)$ is presented by generators and the Weyl relations
$[v,u]=\omega(u,v)$, and $g$ preserves $\omega$, the elements $g(u)$ satisfy
the same relations as the $u$. Hence the substitution $u\mapsto g(u)$ extends
to a unique $k$-algebra endomorphism $\Phi_{g}$ of $A_n(k)$ \cite[Th. 1.3.1, Ex. 7.6.9]{coutinho}. Because $g$ is an automorphism, so is $\Phi_g$, and we have $\Phi_g^{-1} = \Phi_{g^{-1}}.$

\formalref{\leandecl{Stafford38/Paper/Section3/PhiG.lean}{68}{Stafford38.Paper.Section3.phiG},
\leandecl{Stafford38/Paper/Section3/PhiG.lean}{95}{Stafford38.Paper.Section3.step_phi_g}}

\begin{proposition}[{cf. \cite[Thm. 1.1]{stafford1977}}]\label{prop:monic}
Let $0\neq d\in A_n(k)$ with $N:=\deg_B d\ge1$. There exists an automorphism $\Phi_g$ of $A_n(k)$ and
$c\in k^{\times}$ such that
\begin{equation}
  \tilde d:=c^{-1}\Phi_g(d)=p^{\,N}+\sum_{j<N}p^{\,j}a_j, 
  \qquad a_j\in C,\quad \deg_B a_j\le N-j\label{eq:monic-d}
\end{equation}
in the Notation \ref{conv:pair}. $\tilde d$ itself has $\deg_B \tilde d = N.$
\end{proposition}

\begin{proof}
This proposition is a Bernstein-filtrered analogue of Theorem 1.1 of \cite{stafford1977}. For our further proof, we must however be able to control the Berstein degree of $\tilde d$, which the original theorem does not provide.
Let $s = \sigma_N^B(d)\in k[x,\xi]$ be the Bernstein principal symbol of $d$, i.e. s is a nonzero
homogeneous polynomial of degree $N$ on $\mathbb V$. A characteristic-zero field is
infinite, so every nonzero polynomial has a nonvanishing value: we choose $v\in \mathbb V$
with $s(v)\neq0$. Because $s$ is homogeneous of degree $N\geq 1$, $s(0)$ must vanish, therefore $v$ is nonzero. 
Let $e_p\in\mathbb V$ be the standard basis vector corresponding to $p$.
By Lemma~\ref{lem:sympcompl} there is an $\omega$-preserving automorphism
$h$ of $\mathbb V$ with $h(e_p)=v$.\\
If $h$ preserves $\omega$, then $h^{\mathsf T}Jh=J$, where $J$ is the Gram matrix
of $\omega$. In particular $h$ is invertible, and inverting gives
$h^{-1}J^{-1}(h^{\mathsf T})^{-1}=J^{-1}$, i.e.\ $J^{-1}=hJ^{-1}h^{\mathsf T}$.
Since $J^2=-\mathbb{I}$, i.e.\ $J^{-1}=-J$, this is $hJh^{\mathsf T}=J$, so
$g:=h^{\mathsf T}$  again preserves $\omega$. \\
Let $\Phi_g$ be the corresponding automorphism of $A_n(k)$, extending $u\mapsto g(u)$ on
generators and let $M$ be the matrix whose rows are the images of the generators, i.e.\
$g(u_i)=\sum_j M_{ij}u_j$. Since $\Phi_g$ sends each generator to a
degree-one element, it maps $B_m$ into $B_m$ for all $m$, and since
Bernstein symbols are multiplicative, a standard monomial of degree $N$ is
mapped to an element of $B_N$ whose class modulo $B_{N-1}$ is obtained by
substituting $u_i\mapsto\sum_j M_{ij}u_j$ in the corresponding commutative
monomial. Hence $\Phi_g(d)\equiv s(Mu)\bmod B_{N-1}$. \\
As $M$ is invertible,
$s(Mu)$ is again a nonzero homogeneous polynomial of degree $N$, so
$\deg_B\Phi_g(d)=N$ and the Bernstein symbol of $\Phi_g(d)$ is the
function $s_g\colon a\mapsto s(Ma)$. Since the rows of $M$ are the coordinate vectors of the $g(u_i)$, the matrix of
$g$ is $M^T$; as $g=h^T$, so the matrix of $h$ is $M$, i.e.\ $Me_j=h(e_j)$. Hence the
value of $s_g$ at $e_p$ is
\[ s_g(e_p) = s(Me_p)=s(h(e_p))=s(v)\neq0 .\]
We now set $c:=s(v)$. Since $c$ is a nonzero scalar, it is a central unit of
$A_n(k)$, and $\tilde d:=c^{-1}\Phi_g(d)$ satisfies $\deg_B\tilde d=N$ and
$\sigma_N^B(\tilde d)(e_p)=1$.
Using the normal form $A_n(k)=\bigoplus_j p^{\,j}C$, we can write
$\tilde d=\sum_{j = 0}^N p^{\,j}a_j$ with $a_j\in C$. The Bernstein filtration is
compatible with this decomposition, so
$N=\deg_B\tilde d=\max_j\,(j+\deg_B a_j)$, and we have $\deg_B a_j\le N-j$ for
all $j$; in particular $a_j=0$ for $j>N$ and $a_N\in C\cap B_0=k$. 
The symbol is
\[
  \sigma_N^B(\tilde d)=\sum_{j+\deg_B a_j=N}\sigma_1^B(p)^{\,j}\,\sigma_{N-j}^B(a_j),
\]
For $j<N$ the factors $\sigma_N^B(a_j)$ must have a positive degree and must therefore be non-constant in a coordinate other than $p$, so they vanish at $e_p$. Because we defined $\tilde d$ to be $c^{-1}\Phi_g(d)$ , evaluating at $e_p$ gives $a_N=\sigma_N^B(\tilde d)(e_p)=c^{-1}\sigma_N^B(\Phi_g(d))(e_p) = c^{-1}s_g(e_p) = c^{-1}s(v)=1$, and therefore
\[ \tilde d=p^{\,N}+\sum_{j<N}p^{\,j}a_j,\qquad a_j\in C,\quad
  \deg_B a_j\le N-j. \]
\end{proof}

\formalref{%
\leandecl{Stafford38/Paper/Section3/Monic.lean}{243}{Stafford38.Paper.Section3.proposition_3_2},
\leandecl{Stafford38/Paper/Section3/MonicForm.lean}{44}{Stafford38.Paper.Section3.IsMonicNormalForm}}

\begin{remark}\label{rem:transport}
Since $\Phi_g$ is a ring automorphism, the identity $1=\tilde dR+x_1^N\tilde dS$ would transport back to
\[ 1=d\bigl(c^{-1}\Phi_g^{-1}(R)\bigr)+F\,d\bigl(c^{-1}\Phi_g^{-1}(S)\bigr),
  \qquad F=\Phi_g^{-1}(x_1)^N.\]
\end{remark}

\formalref{\leandecl{Stafford38/Paper/Section3/Transport.lean}{77}{Stafford38.Paper.Section3.remark_3_3},
\leandecl{Stafford38/Paper/Section3/Transport.lean}{92}{Stafford38.Paper.Section3.remark_3_3_of_Q_eq_zero}}

We will now define the modules
\begin{equation} \mathcal I=\tilde dA_n+x_1^N\tilde dA_n,\qquad Q=A_n/\mathcal I.
 \label{eq:I-Q}\end{equation}

We now use the order filtration on $A_n(k)$ (where coordinates have
degree zero and derivatives have degree one) to characterize $P =
\sigma^{\mathrm{ord}}_N(\tilde d)$, the order principal symbol of $\tilde
d$ and show that $P$ annihilates the graded module $\gr Q$. 

\begin{definition}\label{def:monic-xi1}
We say a polynomial $P=\sum_{i=0}^Nc_i\,\xi_1^i\in
k[\xi_1,\ldots,\xi_n]$, written as a polynomial in $\xi_1$ with
coefficients $c_i \in k[\xi_2,\ldots,\xi_n]$, is \emph{monic of degree $N$ in
$\xi_1$} if $c_N=1$.
\end{definition}
\formalref{\leandecl{Stafford38/Paper/Section3/OrderSymbol.lean}{52}{Stafford38.Paper.Section3.IsMonicInXi1}}

\begin{lemma}\label{lem:order-symbol}
Let $\tilde d=p^N+\sum_{j<N}p^{\,j}a_j$ with $a_j\in C$, $\deg_Ba_j\le N-j$,
as in \eqref{eq:monic-d}, and set $N=\deg_B\tilde d$. Then $\tilde d$ has
order $N$, and its order principal symbol is a polynomial in $\xi$ alone,
\begin{equation} \sigma_N^{\mathrm{ord}}(\tilde d)=P(\xi),\qquad
 P\in k[\xi_1,\ldots,\xi_n],\qquad \deg P=N.
 \label{eq:P}\end{equation}
$P$ is monic of degree $N$ in $\xi_1$ and in particular
\begin{equation} P(1,0,\ldots,0)=1.
 \label{eq:P-e1}\end{equation}
\end{lemma}

\begin{proof}
Since $\deg_B\tilde d=N$, every standard (PBW) monomial $x^\alpha\partial^\beta$
occurring in $\tilde d$ satisfies $|\alpha|+|\beta|\le N$. A monomial
contributing to the order-$N$ part has $|\beta|=N$, hence $|\alpha|=0$: it
is a pure $\partial^\beta$ with scalar coefficient. Thus $\tilde d$ has
order $N$, and its order-$N$ part is a $k$-linear combination of
$\partial^\beta$ with $|\beta|=N$, giving $\sigma^{\mathrm{ord}}_N(\tilde
d)=P(\xi_1, \dots, \xi_n)$, homogeneous of degree $N$ with coefficients in $k$.

The monomial $p^N=\partial_1^N$ occurs in $\tilde d$ with coefficient $1$.
Every other order-$N$ monomial comes from some $p^ja_j$ with $j<N$ and
$a_j\in C$ (hence involving no further factor of $p=\partial_1$), so it
carries a positive power of some transverse $\partial_i$, $i>1$. Such a
monomial contributes only to coefficients $c_i$ with $i<N$, never to
$c_N$; hence $c_N$ receives no contribution besides the coefficient of
$p^N$ itself, giving $c_N=1$ identically, i.e.\ $P$ is monic of degree $N$
in $\xi_1$. Evaluating at $(1,0,\ldots,0)$ then gives $P(1,0,\ldots,0)=1$.
\end{proof}
\formalref{\leandecl{Stafford38/Paper/Section3/OrderSymbol.lean}{211}{Stafford38.Paper.Section3.lemma_3_5}}

\begin{proposition}\label{prop:charP}
Equip $Q$ with the order filtration $Q_m:=qF_m=\operatorname{im}(F_m\to Q)$
induced by the cyclic generator $q = [1]$, the class of $1$ in $Q$ and write $P=\sigma^{\mathrm{ord}}_N(\tilde d)$.  Then
\begin{equation}\label{eq:P-annihilates}
 P\in\Ann_{k[x,\xi]}\gr Q.
\end{equation}
\end{proposition}
\begin{proof}
$Q_\bullet$ is a good filtration because $Q$ is cyclic. Let $\bar q$ be the class of $q$ in $\gr Q$. We know that $q\,\tilde d=[\tilde d]=0$ in $Q$ (because $\tilde d \in \cal I$) and $\tilde d$ has order $N$ (Lemma \ref{lem:order-symbol}), so $q\tilde d \in Q_N$. Therefore $\bar q\cdot P$ , its image in $Q_N/Q_{N-1}$, must also be zero, and hence $P$ annihilates the class $\bar q$ of $q$ in $\gr Q$.
Since $\gr Q$ is generated over $k[x,\xi]$ by $\bar q$, we get $P\in\Ann_{k[x,\xi]}\gr Q$. 
\end{proof}

\section{$x_1$ as a surjective map}\label{sec:euler}
The first goal of this section will be to show that for $Q$ defined as in Eq. \eqref{eq:I-Q}, we have $Q = Qx_1$ (Proposition \ref{prop:x-surjective}). We will then use this result, as well as Gabber's theorem (Theorem \ref{ithm:gabber}) to establish two contradictory length estimations for the graded module $\gr Q/(\gr Q)x_1$ that we have assumed to be nonzero, concluding that therefore, also $\gr Q = (\gr Q)x_1$ (Corollary \ref{cor:Uzero}).\\

\formalref{\leandecl{Stafford38/Paper/Section3/Transport.lean}{35}{Stafford38.Paper.Section3.idealI}
$=\operatorname{span}\{d,x^Nd\}$ as a right ideal,
\leandecl{Stafford38/Paper/Section3/Transport.lean}{40}{Stafford38.Paper.Section3.Q}}

We define $E=x_1p$. The commutator of $E$ defines a grading of $A_n$: 
$x_1$ has Euler degree $1$, $p$ has Euler degree $-1$, and the transverse Weyl algebra
has Euler degree zero. Let $A_{\geq0}$ denote the nonnegative part of $A_n$ in this grading. In PBW form (with $x_i$ on the left of $\partial_i$), $A_{\geq0}$ is generated by the transverse algebra, $E$, and $x_1$.

By the commutation relations of the Weyl algebra, it holds that
\begin{align}
    Ex_1=x_1(E+1),\qquad x_1E=(E-1)x_1,\nonumber \\
    Ep=p(E-1), \qquad pE=(E+1)p
    \label{eq:nonnegative-pbw}
\end{align}
Additionally, $x_1$ commutes with the transverse algebra, so $x_1$ is normal in $A_{\geq 0}$ (i.e. for every $r \in A_{\geq 0}$, there exists $r'\in A_{\geq 0}$ such that $x_1r=r'x_1$ \cite[4.1.10]{MR}) . This also means that every element of Euler degree $\ge1$ lies in $A_{\ge0}x_1$. 
\formalref{%
\leandecl{Stafford38/Paper/Section4/Euler.lean}{65}{Stafford38.Paper.Section4.eulerDegree},
\leandecl{Stafford38/Paper/Section4/Euler.lean}{111}{Stafford38.Paper.Section4.nonnegativePart},
\leandecl{Stafford38/Paper/Section4/Euler.lean}{129}{Stafford38.Paper.Section4.positivePart},
\leandecl{Stafford38/Paper/Section4/Euler.lean}{637}{Stafford38.Paper.Section4.step_nonnegative_pbw},
\leandecl{Stafford38/Paper/Section4/Euler.lean}{210}{Stafford38.Paper.Section4.step_coordinate_normal}}

\begin{lemma}\label{lem:euler}
For every $m\ge1$,
\begin{equation}\label{eq:eulerprod}
  x_1^{m}p^{\,m}=\prod_{i=0}^{m-1}(E-i), \qquad p^{\,m} x_1^{m}=\prod_{i=1}^{m}(E+i),
\end{equation}
the factors in each product commuting with one another.
\end{lemma}

\begin{proof}
Since both products in \eqref{eq:eulerprod} are polynomials in the single element $E$, their factors commute.
We induct on $m$; both identities are trivial for $m=1$. Assuming the
first identity for $m$,
\[\begin{aligned}
x_1^{m+1}p^{\,m+1}
 &=x_1\Bigl(\prod_{i=0}^{m-1}(E-i)\Bigr)p
 =\Bigl(\prod_{i=0}^{m-1}(E-1-i)\Bigr)x_1p\\
 &=\Bigl(\prod_{i=1}^{m}(E-i)\Bigr)E
 =\prod_{i=0}^{m}(E-i),
\end{aligned}\]
using $x_1E=(E-1)x_1$ repeatedly and the definition of $E$. The second identity follows
identically from $Ep=p(E-1)$.
\end{proof}
\formalref{\leandecl{Stafford38/Paper/Section4/Euler.lean}{348}{Stafford38.Paper.Section4.lemma_4_1}}

\begin{lemma}\label{lem:euler-residue}
There exists an element $u_0\in A_{\geq0}$ such that
\begin{equation} 1+u_0x_1\in \mathcal I.
 \label{eq:euler-residue}\end{equation}
\end{lemma}

\begin{proof}
From Lemma \ref{lem:euler} we define
\[p^Nx_1^N=\prod_{r=1}^N(E+r)=:R_+(E),\qquad
 x_1^Np^N=\prod_{r=0}^{N-1}(E-r)=:R_-(E).\]
The roots of $R_+$ and $R_-$ are disjoint (the sets $\{-1,\ldots,-N\}$ and $\{0,\ldots,N-1\}$ of integers, embedded in the prime field
$\mathbb Q\subset k$), so $\gcd(R_+,R_-)=1$ in $\mathbb Q[E]$. By the
extended Euclidean algorithm, there exist $u,v\in\mathbb Q[E]\subset
k[E]$ with
\[ R_+(E)u(E)+R_-(E)v(E)=p^Nx_1^Nu(E)+x_1^Np^Nv(E)=1.\]
We write $\tilde d=p^N+\Sigma$ with $\Sigma = \sum_{j<N}p^{\,j}a_j$. Then
\[(\tilde d\,x_1^N)u(E)+(x_1^N\tilde d)v(E)=1+w\in \mathcal I,\]
where
\[w=\Sigma x_1^Nu(E)+x_1^N\Sigma v(E).\]
We can write each $a_j$ in $\Sigma$ in powers of $x_1$, with coefficients in the transverse Weyl algebra. A term $p^jx_1^lb$ coming from $p^ja_j$ has $l\geq0$ and $j<N$, and multiplication by $x_1^N$ gives that term Euler degree at least $N-j>0$. $u$ and $v$ only depend on $E$, so they have Euler degree 0. Therefore every term of $w$ has positive Euler degree, and thus $w\in A_{\geq0}x_1$ - or equivalently, there exists a $u_0\in A_{\geq0}$, such that $w=u_0x_1$.
\end{proof}
\formalref{%
\leandecl{Stafford38/Paper/Section4/Residue.lean}{242}{Stafford38.Paper.Section4.lemma_4_2}}

\begin{proposition}\label{prop:x-surjective}
Right multiplication by $x_1$ on $Q$ is surjective.
\end{proposition}
\begin{proof}
Let $q$ be the class of $1$ in $Q$. By Equation \eqref{eq:euler-residue}, we have
$q=-qu_0x_1$, so $q\in(qA_{\geq0})x_1$. We use the normality of $x_1$ to show
\begin{equation} qA_{\geq0}=(qA_{\geq0})x_1.
 \label{eq:normal-surj}\end{equation}
The inclusion $\supseteq$ holds because $x_1\in A_{\geq0}$ and
$A_{\geq0}$ is a ring. For $\subseteq$, let $r\in A_{\geq0}$ and choose
$r'$ with $x_1r=r'x_1$ (possible since $x_1$ is normal) then $qr=-(qu_0)x_1r=-(qu_0r')x_1\in(qA_{\geq0})x_1$ since $u_0r'\in A_{\geq0}$.

$qA_{\geq0}$ is stable under right multiplication by $p$: If $m=m'x_1$, then
\[ mp=m'(x_1p)=m'E\in qA_{\geq0}.\]
It is of course also stable under $C$, since $C=A_{n-1}(k)[x_1]\subseteq A_{\ge0}$ because no element of $C$ contains $p$. Since $C$ and $p$ together generate $A_n$, we have $qA_{\geq0}=qA_n=Q$. Equation \eqref{eq:normal-surj} is therefore $Q=Qx_1$.
\end{proof}
\formalref{\leandecl{Stafford38/Paper/Section4/Surjectivity.lean}{223}{Stafford38.Paper.Section4.proposition_4_3}}

\begin{definition}[{\cite{Stacks}[Tag 00CM]}] \label{def:localization}
Let $R$ be a commutative ring and $S\subseteq R$ a multiplicative subset
(i.e.\ $1\in S$ and $S$ is closed under products). The \emph{localization}
$S^{-1}R$ is the ring of fractions $a/s$ ($a\in R$, $s\in S$), where
\[ \frac{a}{s}=\frac{a'}{s'}\iff t(s'a-sa')=0\ \text{ for some }t\in S,\]
with the usual addition and multiplication of fractions. Analogously, for an $R$-module
$N$, the \emph{localization} $S^{-1}N$ is the $S^{-1}R$-module of fractions
$y/s$ ($y\in N$, $s\in S$), with $y/s=y'/s'$ iff $t(s'y-sy')=0$ for some
$t\in S$. For a prime $\mathfrak p\subset R$ we write
$R_{\mathfrak p}=S^{-1}R$ and $N_{\mathfrak p}=S^{-1}N$ for
$S=R\setminus\mathfrak p$. Localization is exact: it sends exact sequences
of $R$-modules to exact sequences of $S^{-1}R$-modules \cite[Tag 00CS]{Stacks}.
\end{definition}
\begin{definition}[{\cite{Stacks}[Definition 10.40.1]}]\label{def:supp}
    Let $R$ be a ring and let $M$ be an $R$-module. The \emph{support} of $M$ is the set 
    \[\Supp_R(M) = \{\mathfrak{p} \in \Spec(R) | M_{\mathfrak{p}} \neq 0\}.\]
\end{definition}

\begin{definition}[\mbox{\cite{coutinho, cohn}, \cite[Tags~02LY, 00J3, 00IV]{Stacks}}]\label{def:length}
    The length $\ell$ of a $T$-module of finite length $N$ is defined as the number of steps in a composition series of $N$ such that
    \[0 = N_0\subsetneq N_1\subsetneq \dots \subsetneq N_{\ell} = N\]
    where each successive quotient $N_i/N_{i-1}$ is simple (i.e. has no proper nontrivial submodules). The length is well-defined, as it does not depend on the choice of the composition series \cite[Tag 00J3]{Stacks}. Lengths are additive for short exact sequences.
\end{definition}
\formalref{\leandecl{Stafford38/Paper/Section4/Support.lean}{59}{Stafford38.Paper.Section4.definition_4_6};
\texttt{Module.length} (Mathlib); finite length of the localized modules:
\leanlib{AlgebraicAnalysis/Module/MinimalSupportKernelCokernelLengths.lean}{48}{AlgebraicAnalysis.localized_kernel_and_cokernel_isFiniteLength}}
\begin{lemma}\label{lem:supp-criterion}
Let $R$ be a commutative ring, $N$ an $R$-module and $\mathfrak p\subset R$ a
prime. Then $\mathfrak p\in\operatorname{Supp}_R(N)$ if and only if there is
$y\in N$ with $\operatorname{Ann}_R(y)\subseteq\mathfrak p$.
\end{lemma}
\formalref{\leandecl{Stafford38/Paper/Section4/Support.lean}{111}{Stafford38.Paper.Section4.lemma_4_7}}
\begin{proof}
In $N_{\mathfrak p}$, a
fraction is zero when it can be cleared by an element outside $\mathfrak p$.
For $y/1$ this means: there is $u\in R\setminus\mathfrak p$ with $uy=0$.
Such a $u$ lies in $\operatorname{Ann}_R(y)$ and not in $\mathfrak p$. So
\[y/1=0\iff\operatorname{Ann}_R(y)\not\subseteq\mathfrak p,\]
and, taking the contrapositive,
\[y/1\neq0\iff\operatorname{Ann}_R(y)\subseteq\mathfrak p.\]
\emph{$(\Leftarrow)$} Suppose $y\in N$ has
$\operatorname{Ann}_R(y)\subseteq\mathfrak p$. Then, $y/1\neq0$ in
$N_{\mathfrak p}$. So $N_{\mathfrak p}$ has a nonzero element, hence
$N_{\mathfrak p}\neq0$ and $\mathfrak p\in\operatorname{Supp}_R(N)$.\\
\emph{$(\Rightarrow)$} Suppose $N_{\mathfrak p}\neq0$. Pick a nonzero
element of $N_{\mathfrak p}$. It has the form $y/s$ with $y\in N$ and
$s\notin\mathfrak p$. The element $s$ is invertible in $N_{\mathfrak p}$, so
\[y/s=0\iff y/1=0.\]
Since $y/s\neq0$, we get $y/1\neq0$, and therefore
$\operatorname{Ann}_R(y)\subseteq\mathfrak p$, which is the desired $y$.
\end{proof}

For the next part, we write $M=\operatorname{gr}Q$ and denote the right multiplication by $x_1$ by $f:M\rightarrow M$. We furthermore define
\[
 \begin{gathered}
 V=\ker(f),\qquad U=\operatorname{coker}(f) =M/x_1M\\
 T=k[x_2,\ldots,x_n,\xi_2,\ldots,\xi_n],
 \end{gathered}
\]
hence $T[x_1,\xi_1] = k[x,\xi]$. If $n=1$, we define $T=k$.
We want to prove that $U=0$ by comparing
two lengths over a localization of $T$.

\begin{lemma}\label{lem:tangential-finiteness}
The modules $U$ and $V$ are finite over $T$, and
$\operatorname{Supp}_T V\subseteq\operatorname{Supp}_T U$.
If $\mathfrak q$ is a minimal prime in $\operatorname{Supp}_T U$, then
$U_{\mathfrak q}$ and $V_{\mathfrak q}$ have finite length over
$T_{\mathfrak q}$.
\end{lemma}

\begin{proof}
The good quotient filtration makes $M = \operatorname{gr} Q$ finitely generated over $\operatorname{gr} A_n = k[x,\xi] = T[{x_1}, \xi_1]$.
Since $P$ annihilates $M$ (Proposition~\ref{prop:charP}), the action of
$T[x_1,\xi_1]$ on $M$ factors through $T[x_1,\xi_1]/(P)$. As $P$ is monic of
degree $N$ in $\xi_1$, in $T[x_1,\xi_1]/(P)$ we have $\xi_1^N=-\sum_{i<N}c_i\xi_1^i$ with
$c_i\in T$, so $T[x_1,\xi_1]/(P)$ is generated as a $T[x_1]$-module by $1,\xi_1,\dots,\xi_1^{N-1}$.
If $m_1,\dots,m_l$ generate $M$ over $T[x_1,\xi_1]/(P)$, then the products $\xi_1^im_s$
($0\le i<N$, $1\le s\le l$) generate $M$ over $T[x_1]$. Hence $M$ is finite over
$T[x_1]$.
Because $M$ is finite, and $T[{x_1}]$ is a Noetherian ring, the kernel $V$ and cokernel $U$ of multiplication in $M$ by ${x_1}$ are also finite over $T[{x_1}]$. Since they are annihilated by ${x_1}$, they are finite over $T$. 

Let $\mathfrak p\in\operatorname{Supp}_T(V)$. By Lemma \ref{lem:supp-criterion}, we can choose $v\in V=\ker f$ with
$\operatorname{Ann}_T(v)\subseteq\mathfrak p$. Let us now consider the annihilators of $v$ in $k[x,\xi]=T[x_1, \xi_1]$, defining $S$ as the image of $T\setminus \mathfrak p$ in $k[x,\xi]$: $\operatorname{Ann}_{k[x,\xi]}(v)\cap S = \varnothing$, because were $r \in T\setminus\mathfrak p$ and in $\operatorname{Ann}_{k[x,\xi]}(v)$, then $v$ would be annihilated by $r$, which in turn would imply that $r$ is also in $\operatorname{Ann}_T(v)$. Since $\operatorname{Ann}_T(v) \subseteq \mathfrak p$, but we assumed $r \in T\setminus\mathfrak p$, this is a contradiction. There must therefore exist a prime $\mathfrak Q \in k[x,\xi]$ with $\mathfrak Q \cap S = \varnothing$ and $\operatorname{Ann}_{k[x,\xi]}(v) \subseteq \mathfrak Q$. Since $v \in V = \ker f$, by Lemma \ref{lem:supp-criterion}, it follows that $\mathfrak Q \in \Supp_{k[x,\xi]}(\ker f).$\\
From Vasconcelos' theorem \cite[Theorem~2.4]{matsumura} \cite{vasconcelos} follows that $\operatorname{Supp}_{k[x,\xi]}(\ker f)\subseteq \operatorname{Supp}_{k[x,\xi]}(\operatorname{coker}f)$: 
We prove the contrapositive. Let $\mathfrak Q$ be a prime of ${k[x,\xi]}$ with
$\mathfrak Q\notin\operatorname{Supp}_{k[x,\xi]}(\operatorname{coker}f)$, i.e.\
$(\operatorname{coker}f)_{\mathfrak Q}=0$.
Since localization is exact, $(\operatorname{coker}f)_{\mathfrak Q}=
\operatorname{coker}(f_{\mathfrak Q})$, so
$\operatorname{coker}(f_{\mathfrak Q})=0$: the endomorphism
$f_{\mathfrak Q}\colon M_{\mathfrak Q}\to M_{\mathfrak Q}$ is therefore surjective.
The localization $M_{\mathfrak Q}$ is a finite module over the Noetherian ring
${k[x,\xi]}_{\mathfrak Q}$, and $f_{\mathfrak Q}$ is ${k[x,\xi]}_{\mathfrak Q}$-linear, because
$f$ is multiplication by $x_1\in {k[x,\xi]}$ and ${k[x,\xi]}$ is commutative. By Vasconcelos'
theorem, a surjective endomorphism of a finite module is injective. Hence
$f_{\mathfrak Q}$ is injective, and by exactness of localization,
$(\ker f)_{\mathfrak Q}=\ker(f_{\mathfrak Q})=0$, i.e.\
$\mathfrak Q\notin\operatorname{Supp}_{k[x,\xi]}(\ker f)$. The contrapositive is therefore $\mathfrak Q\in\operatorname{Supp}_{k[x,\xi]}(\ker f)$,
then $\mathfrak Q\in\operatorname{Supp}_{k[x,\xi]}(\operatorname{coker}f) = \Supp_{k[x,\xi]}(U).$ By
Lemma~\ref{lem:supp-criterion} there therefore exists an $u\in U$ with
$\operatorname{Ann}_{k[x,\xi]}(u)\subseteq\mathfrak Q$. 
To get back to the support over $T$, we define $r\in\operatorname{Ann}_T(u)$. Then $ru=0$, so
$r\in\operatorname{Ann}_{k[x,\xi]}(u)\subseteq\mathfrak Q$. If $r\notin\mathfrak p$,
the image of $r$ would lie in $\mathfrak Q\cap S=\varnothing$, which is
impossible. Hence $r\in\mathfrak p$, i.e.\
$\operatorname{Ann}_T(u)\subseteq\mathfrak p$, and
Lemma~\ref{lem:supp-criterion} gives
$\mathfrak p\in\operatorname{Supp}_T(U)$. Since $\mathfrak p\in
\operatorname{Supp}_T(V)$ was arbitrary,
$\operatorname{Supp}_T(V)\subseteq\operatorname{Supp}_T(U)$.

Let now $\mathfrak q$ be a minimal prime in $\operatorname{Supp}_TU$. The primes of
$T_{\mathfrak q}$ correspond to the primes $\mathfrak p'\subseteq\mathfrak q$ of $T$\cite[Prop 3.11]{atiyah}, and
$(U_{\mathfrak q})_{\mathfrak p'}=U_{\mathfrak p'}$. Hence $\operatorname{Supp}(U_{\mathfrak q})$
corresponds to the primes $\mathfrak p'\subseteq\mathfrak q$ lying in $\operatorname{Supp}_TU$;
by minimality of $\mathfrak q$ the only such prime is $\mathfrak q$ itself, so
$\operatorname{Supp}(U_{\mathfrak q})=\{\mathfrak qT_{\mathfrak q}\}$, the set consisting
only of the maximal ideal of $T_{\mathfrak q}$. The same argument applies to $V$: a prime
$\mathfrak p'\subseteq\mathfrak q$ in $\operatorname{Supp}_TV$ lies in $\operatorname{Supp}_TU$
by the inclusion just proved, hence equals $\mathfrak q$, so
$\operatorname{Supp}(V_{\mathfrak q})\subseteq\{\mathfrak qT_{\mathfrak q}\}$.

Now $T_{\mathfrak q}$ is a Noetherian local ring with maximal ideal $\mathfrak qT_{\mathfrak q}$,
and $U_{\mathfrak q}$, $V_{\mathfrak q}$ are finite modules (localizations of the finite
$T$-modules $U$, $V$) whose support is contained in $\{\mathfrak qT_{\mathfrak q}\}$. For a
finite module $L$ over $T_{\mathfrak q}$ $\operatorname{Supp}L$ is equal to the vanishing locus of its annihilators \cite[Tag 00L2]{Stacks},
so the only prime containing $\operatorname{Ann}L$ is the maximal ideal (if $L\neq0$). Hence
$T_{\mathfrak q}/\operatorname{Ann}L$ is a Noetherian ring in which every prime is maximal,
i.e.\ an Artinian ring, and a finite module over an Artinian ring has finite length \cite[Thm. 8.5]{atiyah}\cite[Tag 00KH]{Stacks}. Thus
$U_{\mathfrak q}$ and $V_{\mathfrak q}$ have finite length over $T_{\mathfrak q}$
($V_{\mathfrak q}$ possibly zero).
\end{proof}
\formalref{\leandecl{Stafford38/Paper/Section4/TangentialFiniteness.lean}{164}{Stafford38.Paper.Section4.lemma_4_8}}

\begin{lemma}\label{lem:page-inequality}
For every minimal prime $\mathfrak q$ in $\operatorname{Supp}_T U$,
\[
 \ell(U_{\mathfrak q})
 \le \ell(V_{\mathfrak q})
\] over $T_{\mathfrak q}.$
\end{lemma}

\begin{proof}
To prove the inequality we set up a spectral sequence measuring how $\ker(f)$ and
$\operatorname{coker}(f)$ acting on $Q$ itself compare, order by order, to $V$ and
$U$, acting on $\gr Q$.

Consider the filtered complex $Q\xrightarrow{f}Q$, with decreasing filtration $G^p=F_{-p}Q$, where $F$ is
the order filtration. Since ${x_1}$ has order $0$, $f(G^p)\subseteq G^p$ for every
$p$; the induced map on $G^p/G^{p+1}$, summed over $p$, is exactly the action of
${x_1}$ on $M=\operatorname{gr}Q$. This is the crudest approximation to $\ker_Q(f)$ and
$\operatorname{coker}_Q(f)$: it correctly identifies $V,U$ as approximations, but it
does not yet see elements that vanish to higher order than this naive count
predicts. The spectral sequence refines this approximation in successive stages,
called \emph{pages} (see also \cite[Construction 5.4.6]{weibel}): for each $r\ge0$, the source and target pages are
\[
 A_r=\bigoplus_{p\in\mathbb Z}\frac{G^p\cap f^{-1}(G^{p+r})}{G^{p+1}\cap f^{-1}(G^{p+r})},
 \qquad
 B_r=\bigoplus_{p\in\mathbb Z}\frac{G^p}{(G^p\cap f(G^{p-r+1}))+G^{p+1}},
\]
where the degree-$p$ summand of $A_r$ isolates elements of $G^p$ whose image under
$f$ drops $r$ extra filtration levels (candidates for the kernel, accurate to
order $r$), and the degree-$p$ summand of $B_r$ isolates the part of $G^p$ not yet
hit by $f$ from $r-1$ levels shallower.

Applying $f$ to a representative gives a well-defined map $d_r:A_r\to B_r$,
shifting the summand index by $r$; one checks directly from these definitions that
\[ A_0=B_0=M,\qquad d_0={x_1},\qquad A_{r+1}=\ker d_r,\qquad B_{r+1}=\operatorname{coker}d_r,\]
so that each page is obtained from the previous one by taking kernel and cokernel
of $d_r$, exactly as in the construction of a spectral sequence from a filtered
complex.

These pages carry a $T$-module structure. The pages $A_r,B_r$ are
subquotients of $Q$, which a priori carries an action of the transverse Weyl
algebra $A_{n-1}(k)$ by right multiplication. This action commutes with
$f$, since the transverse algebra commutes with $x_1$, so the maps
$d_r$ are compatible with it. Two transverse elements commute only up to a
commutator, which lowers the order filtration by one. On the page
quotients such lower-order terms vanish, so the induced actions commute, and
the commutative ring $T$ acts on every $A_r$ and $B_r$, via the symbols of
the transverse generators. The maps $d_r$ and the successor-page maps are
$T$-linear. In particular, this means that localizing at $T\setminus\mathfrak q$ and
measuring lengths over $T_{\mathfrak q}$ is a well-defined action for the pages.

Now we fix a minimal prime $\mathfrak q$ in $\operatorname{Supp}_TU$ and localize at
$T\setminus\mathfrak q$. We localize the $T$-modules
$A_r,B_r$ and their $T$-linear maps. Exactness of localization preserves the
successor-page kernel and cokernel identities. By
Lemma~\ref{lem:tangential-finiteness}, $(A_1)_{\mathfrak q}=V_{\mathfrak q}$ and
$(B_1)_{\mathfrak q}=U_{\mathfrak q}$ have finite length over $T_{\mathfrak q}$; since each
$(A_{r+1})_{\mathfrak q}\subseteq (A_r)_{\mathfrak q}$ and $(B_{r+1})_{\mathfrak q}$ is a quotient of $(B_r)_{\mathfrak q}$, all later pages
have finite length as well.

The maps $(B_1)_{\mathfrak q}\to (B_r)_{\mathfrak q}$ are surjective by construction.
Let $K_r:=\ker((B_1)_{\mathfrak q}\to (B_r)_{\mathfrak q})$. Since $(B_1)_{\mathfrak q}\to (B_{r+1})_{\mathfrak q}$ factors through $(B_1)_{\mathfrak q}\to (B_r)_{\mathfrak q}\to (B_{r+1})_{\mathfrak q}$, we have
$K_r\subseteq K_{r+1}$, so $(K_r)_{r\ge1}$ is an increasing chain of submodules of $(B_1)_{\mathfrak q}$.
This chain exhausts $(B_1)_{\mathfrak q}$: Before localization, take a homogeneous
class in $B_1$, represented by $z\in G^p$. Since $Q=Qx_1$ by Proposition \ref{prop:x-surjective}, choose one preimage
$z=wx_1$. Exhaustiveness of the filtration means that $w$ must be in $G^s$ for some $s$; choose $r\ge1$ with
$p-r+1\le s$. Then $w\in G^{p-r+1}$, so $z$ belongs to the denominator of
$B_r$ and its class vanishes. Every class is a finite sum of homogeneous
classes and hence vanishes on a sufficiently late page. The same holds for a
localized class $z/a$, since its numerator already vanishes on such a page.
As $z\in B_1$ was arbitrary, so is $z/a$ in $(B_1)_{\mathfrak q}$, and $\bigcup_{r\ge1}K_r=(B_1)_{\mathfrak q}$.
Since $(B_1)_{\mathfrak q}$ has finite length, it is a Noetherian $T_{\mathfrak q}$-module. By the ascending
chain condition, the increasing chain $K_1\subseteq K_2\subseteq\cdots$ stabilizes: there exists a finite $R$ with
$K_R=K_{R+i}$ for all $i\geq 0$. Combined with $\bigcup_rK_r=(B_1)_{\mathfrak q}$, stabilization forces $K_R=(B_1)_{\mathfrak q}$.
Finally, $K_R=(B_1)_{\mathfrak q}$ means $(B_1)_{\mathfrak q}\to (B_R)_{\mathfrak q}$ must be the zero map: by the first isomorphism theorem,
$(B_1)_{\mathfrak q}/K_R\cong\operatorname{im}((B_1)_{\mathfrak q}\to (B_R)_{\mathfrak q})$, and $K_R=(B_1)_{\mathfrak q}$ gives $(B_1)_{\mathfrak q}/K_R=0$, so
$\operatorname{im}((B_1)_{\mathfrak q}\to (B_R)_{\mathfrak q})=0$. But $(B_1)_{\mathfrak q}\to (B_R)_{\mathfrak q}$ is also surjective, so $(B_R)_{\mathfrak q}=0$.\\

At each stage $r\ge1$ the differential $(d_r)_{\mathfrak q}$ gives an exact sequence
\[ 0\to (A_{r+1})_{\mathfrak q}\to (A_{r})_{\mathfrak q}\xrightarrow{(d_r)_{\mathfrak q}}(B_{r})_{\mathfrak q}\to (B_{r+1})_{\mathfrak q}\to0,\]
which splits at $\operatorname{im}((d_r)_{\mathfrak q})$ into two short exact sequences.
Additivity of length gives
$\ell((A_{r})_{\mathfrak q})=\ell((A_{r+1})_{\mathfrak q})+\ell(\operatorname{im}(d_r)_{\mathfrak q})$ and
$\ell((B_{r})_{\mathfrak q})=\ell(\operatorname{im}(d_r)_{\mathfrak q})+\ell((B_{r+1})_{\mathfrak q})$; eliminating
$\ell(\operatorname{im}(d_r)_{\mathfrak q})$ yields, for every $r\ge1$,
\begin{equation}\label{eq:additivity}
 \ell((B_r)_{\mathfrak q})+\ell((A_{r+1})_{\mathfrak q})=\ell((A_{r})_{\mathfrak q})+\ell((B_{r+1})_{\mathfrak q}).
\end{equation}

Summing this identity over $r=1,\ldots,R-1$ and reindexing the sums, most terms cancel on both sides, leaving 
\[\ell((B_1)_{\mathfrak q})+\ell((A_R)_{\mathfrak q})=\ell((A_{1})_{\mathfrak q})+\ell((B_{R})_{\mathfrak q}).\]
Since $(B_{R})_{\mathfrak q}=0$ and $(A_{1})_{\mathfrak q}=V_{\mathfrak q}$, $(B_{1})_{\mathfrak q}=U_{\mathfrak q}$, this reads
\[ \ell(U_{\mathfrak q})+\ell((A_{R})_{\mathfrak q})=\ell(V_{\mathfrak q}).\]
As $\ell((A_{R})_{\mathfrak q})\ge0$, we conclude
\[ \ell(U_{\mathfrak q})\le\ell(V_{\mathfrak q}).\]
\end{proof}

\formalref{\leandecl{Stafford38/Paper/Section4/PageInequality.lean}{173}{Stafford38.Paper.Section4.lemma_4_9};
\leanlib{AlgebraicAnalysis/Module/UniformBoundaryVanishing.lean}{106}{AlgebraicAnalysis.exists_uniform_subsingleton_localized}}

\noindent On the other hand, we can also establish a different inequality:

\begin{lemma}\label{lem:strict-support-inequality}
For every minimal prime $\mathfrak q$ in $\operatorname{Supp}_T U$,
\[
 \ell(U_{\mathfrak q})
 > \ell(V_{\mathfrak q})
\] over $T_{\mathfrak q}$.
\end{lemma}
\begin{proof}
We start by showing, that by Gabber's involutivity theorem, $x_1$ cannot lie in any minimal prime. We then define a module $K_{\infty}$ as the kernel of multiplication by a power of $x_1$ over the localized $M_\mathfrak q$, which allows us to split the length of $U_{\mathfrak q}$ into two parts, $\ell(U_{\mathfrak q}) = \ell(V_{\mathfrak q}) + \ell\bigl((M_{\mathfrak q}/K_{\infty})/({x_1}(M_{\mathfrak q}/K_{\infty}))\bigr)$. It remains to show that the second term is nonzero. For this we use the fact that $x_1$ does not lie in any minimal prime to deduce that $x_1$ must therefore be invertible at the minimal primes while lying in every prime of the support of $U_{\mathfrak q}$. Nakayama's lemma then shows that  $(M_\mathfrak q/K_{\infty})/x_1(M_\mathfrak q/K_{\infty})\neq0$. \\

Let us consider $\operatorname{Ann}_{k[x,\xi]} M$. By Gabber's Theorem (\ref{ithm:gabber}, \cite{Gabber1981}), $J = \sqrt{\operatorname{Ann}_{k[x,\xi]} M}$ is involutive.
If $\mathfrak p$ is a minimal prime of $J$,
then $\mathfrak p$ is itself involutive, i.e. $\{\mathfrak p,\mathfrak p\}\subseteq\mathfrak p$ \cite[Thm. D.3.4]{HTT}.
\formalref{%
\leandecl{Stafford38/Paper/Section4/StrictInequality.lean}{405}{Stafford38.Paper.Section4.lemma_4_10}}
From the involutivity of $\mathfrak p$ follows that no $\mathfrak p$ contains ${x_1}$: Assume ${x_1}\in \mathfrak p$. Since $P \in \mathfrak p$,  also the Poisson brackets
\[ \{{x_1}, P\}, \quad \{{x_1}, \{x_1, P\}\}, \dots\]
lie in $\mathfrak p$. $\{{x_1},-\}$ differentiates once with respect to $\xi_1$, and because $P$ is homogeneous of degree $N$ and monic in $\xi_1$, its $N$th derivative is the nonzero constant $N!$ - therefore the $N$th bracket would put $N!$ in $\mathfrak p$. Since $\mathfrak p$ is a proper prime ideal, no nonzero constant can be in $\mathfrak p$, so we must conclude ${x_1} \notin \mathfrak p$.
This persists after localizing at $S=T\setminus\mathfrak q$: write $k[x,\xi]_{\mathfrak q}:=S^{-1}k[x,\xi]$. Since $M$ is finite, $\operatorname{Ann}(M_{\mathfrak q})=S^{-1}\operatorname{Ann}(M)$ \cite[Prop 3.14]{atiyah}, and primes of $k[x,\xi]_{\mathfrak q}$ correspond to primes of $k[x,\xi]$ disjoint from $S$ \cite[Prop 3.11(iv)]{atiyah}. A prime minimal among those containing $J$ and disjoint from $S$ is minimal over $J$ (any smaller prime containing $J$ is again disjoint from $S$). Hence every minimal prime of $\operatorname{Supp}(M_{\mathfrak q})$ is of the form $S^{-1}\mathfrak p$ with $\mathfrak p$ a minimal prime of $J$, and does not contain ${x_1}$.

Consider the increasing chain of kernels of powers of the localization $f_{\mathfrak q}$ of right multiplication with $x_1$
$\operatorname{ker}(f_{\mathfrak q})$ $\subseteq \operatorname{ker}(f^2_{\mathfrak q})$ $\subseteq \dots$.
Since $M_{\mathfrak q}$ is Noetherian, this chain stabilizes. Let
$K_{\infty}:=\ker(f^m_{\mathfrak q})$ at a stage where it stabilizes, so $f_{\mathfrak q}(K_{\infty}) = {x_1}K_{\infty}\subseteq K_{\infty}$ and
$\ker(f^{m+1}_{\mathfrak q})=\ker(f_{\mathfrak q}^m)=K_{\infty}$. Each quotient $\ker(f_{\mathfrak q}^{i+1})/\ker(f_{\mathfrak q}^i)$ maps to
$V_{\mathfrak q}=\ker(f_{\mathfrak q})$ injectively via $v\mapsto {x_1}^iv$, so, as $V_{\mathfrak q}$ has
finite length, so does each quotient and hence $K_{\infty}$ itself.

Two elementary facts about $K_{\infty}$:
\begin{itemize}
\item ${f_{\mathfrak q}}^{-1}K_{\infty}=K_{\infty}$: if $v\in K_{\infty}$ then $f_{\mathfrak q}(v)\in K_{\infty}$ (as $f_{\mathfrak q}(K_{\infty})\subseteq K_{\infty}$); if $f_{\mathfrak q}(v)\in K_{\infty}=\ker(f_{\mathfrak q}^m)$
then $f_{\mathfrak q}^{m+1}(v)=0$, so $v\in\ker({f_{\mathfrak q}}^{m+1})=K_{\infty}$.
\item $K_{\infty}\cap {x_1}M_{\mathfrak q}={x_1}K_{\infty} = f_{\mathfrak q}(K_{\infty})$: "$\supseteq$" is trivial. For "$\subseteq$", if $w=f_{\mathfrak q}(v)\in K_{\infty}$ for
some $v\in M_{\mathfrak q}$, then $v\in f_{\mathfrak q}^{-1}K_{\infty}=K_{\infty}$, so $w\in f_{\mathfrak q}(K_{\infty}) = x_1K_{\infty}$.
\end{itemize}

For any endomorphism $\varphi$ of a finite-length module $L$,
additivity along $0\to\ker \varphi\to L\to \varphi(L)\to0$ and $0\to \varphi(L)\to L\to L/\varphi(L)\to0$ gives
$\ell(\ker \varphi)=\ell(L/\varphi(L))$. We can apply this to $\varphi=f_{\mathfrak q}|_{K_{\infty}}$: since $\ker(f_{\mathfrak q})\subseteq K_{\infty}$, we
have $\ker(f_{\mathfrak q}|_{K_{\infty}})=\ker(f_{\mathfrak q})=V_{\mathfrak q}$, so with $f_{\mathfrak q}(K_{\infty}) = x_1K_{\infty},$
\[\ell(K_{\infty}/{x_1}K_{\infty})=\ell(V_{\mathfrak q}).\]

On the other hand, the natural surjection $M_{\mathfrak q}/({x_1}M_{\mathfrak q})\to (M_{\mathfrak q}/K_{\infty})/({x_1}(M_{\mathfrak q}/K_{\infty}))$ has
the kernel $(K_{\infty}+({x_1}M_{\mathfrak q}))/({x_1}M_{\mathfrak q})\cong K_{\infty}/(K_{\infty}\cap ({x_1}M_{\mathfrak q}))=K_{\infty}/({x_1}K_{\infty})$, giving (we now write $\tilde E = (M_{\mathfrak q}/K_{\infty})$ for better readability) the exact sequence
\[0\longrightarrow K_{\infty}/({x_1}K_{\infty})\longrightarrow U_{\mathfrak q}\longrightarrow
 \tilde E/({x_1}\tilde E)\longrightarrow0,\]
so by additivity and using $\ell(K_{\infty}/({x_1}K_{\infty}))=\ell(V_{\mathfrak q})$,
\[\ell(U_{\mathfrak q}) = \ell(V_{\mathfrak q}) + \ell\bigl(\tilde E/({x_1}\tilde E)\bigr).\]
It remains to show the last term is positive.

Since $\mathfrak q\in\operatorname{Supp}_TU$,
$U_{\mathfrak q}=M_{\mathfrak q}/({x_1}M_{\mathfrak q})\ne0$. We choose $\mathfrak P\in\operatorname{Supp}_{k[x,\xi]_{\mathfrak q}}(M_{\mathfrak q}/({x_1}M_{\mathfrak q}))$
and a prime $\mathfrak p'\subseteq\mathfrak P$ minimal in $\operatorname{Supp}(M_{\mathfrak q})$ \cite[Tag 00E0]{Stacks}. By the above, ${x_1}\notin\mathfrak p'$, so ${x_1}$ is a unit in the localization at $\mathfrak p'$; as $K_{\infty}$ is annihilated by a power of ${x_1}$, $(K_{\infty})_{\mathfrak p'}=0$. From
$0\to (K_{\infty})_{\mathfrak p'}\to (M_{\mathfrak q})_{\mathfrak p'}\to \tilde E_{\mathfrak p'}\to0$ we get
$\tilde E_{\mathfrak p'}\cong (M_{\mathfrak q})_{\mathfrak p'}\ne0$, so $\mathfrak p'$ lies in the support of $M_{\mathfrak q}/K_{\infty}$ in $k[x,\xi]_{\mathfrak q}$.
Since support is closed under passing to larger primes \cite[Tag 00L2]{Stacks} and
$\mathfrak p'\subseteq\mathfrak P$, also $\mathfrak P\in\operatorname{Supp}_{k[x,\xi]_{\mathfrak q}}\tilde E$, i.e.
$\tilde E_{\mathfrak P}\ne0$.

Since $\mathfrak P\in\operatorname{Supp}_{k[x,\xi]_{\mathfrak q}}(M_{\mathfrak q}/({x_1}M_{\mathfrak q}))$ and $M_{\mathfrak q}/({x_1}M_{\mathfrak q})$ is annihilated by
${x_1}$, we have ${x_1}\in\mathfrak P$. If $({x_1}\tilde E_{\mathfrak P})=\tilde E_{\mathfrak P}$, then as
${x_1}$ lies in the maximal ideal of the local ring $(k[x,\xi]_{\mathfrak q})_{\mathfrak P}$, Nakayama's lemma \cite{nakayama}, (see also \cite[Tag~00DV]{Stacks})
would force $\tilde E_{\mathfrak P}=0$, a contradiction. Hence
$\bigl(\tilde E/({x_1}\tilde E)\bigr)_{\mathfrak P}\ne0$, so $\tilde E/({x_1}\tilde E)\ne0$ and $\ell\left(\tilde E/({x_1}\tilde E)\right) > 0$.

Therefore $\ell(U_{\mathfrak q})=\ell(V_{\mathfrak q})+\ell\bigl(\tilde E/{x_1}\tilde E\bigr)
>\ell(V_{\mathfrak q})$.
\end{proof}
\formalref{\leandecl{Stafford38/Paper/Section4/StrictInequality.lean}{405}{Stafford38.Paper.Section4.lemma_4_10}}

\begin{corollary} 
\label{cor:Uzero}
Multiplication by $x_1$ on $\operatorname{gr}Q$ is surjective, i.e.\
$\gr Q=(\gr Q)\,x_1$.
\end{corollary}
\begin{proof}
We write $M = \gr Q$ and $U = M/{x_1}M$. If $U\ne0$, it has a minimal support prime $\mathfrak q$.
Lemma~\ref{lem:tangential-finiteness} makes the relevant lengths finite,
and Lemmas~\ref{lem:page-inequality} and
\ref{lem:strict-support-inequality} contradict one another.
Thus $U = M/x_1M=0$ and therefore $\gr Q=(\gr Q)x_1$.
\end{proof}
\formalref{\leandecl{Stafford38/Paper/Section4/Vanishing.lean}{45}{Stafford38.Paper.Section4.corollary_4_11}}

\section{Avoiding the hyperplane $\{x_1 = 0\}$}\label{sec:conormal-review}
This section uses mainly vocabulary and concepts from algebraic geometry. We will define the most important quantities used here (for a comprehensive overview of the concepts of algebraic geometry, we refer to textbooks such as \cite{hartshorne}), then state the two main results: Theorem~\ref{thm:asymptotic-conormal-formal}
shows that for a closed irreducible subvariety $Y\subseteq X$ not meeting the hyperplane
$\{x_1=0\}$, the covector $dx_1$ is an asymptotic conormal direction of $Y$.
Theorem~\ref{thm:coisotropic-exclusion} uses it to show that no nonempty closed,
fibre-conical, involutive subset of the cotangent space can lie over the complement
of the hyperplane $\{x_1=0\}$ and be contained in the zero set of a polynomial
$P(\xi)$ with $P(1,0,\dots,0)\neq0$.
We prove the first theorem in Subsection~\ref{sub:asym}, with the lemmas stated and
proved where they are used, and the second in Subsection~\ref{sub:excl}.
\subsection{Preliminaries and definitions from algebraic geometry}\label{sub:def_geom}
Throughout this section $K$ denotes an algebraically closed field of characteristic
zero, and we use the notation of Definition~\ref{def:cotangent} with $k=K$. A prime
ideal that is maximal is called a \emph{closed point}. By Theorem~\ref{thm:nss}(a),
the map $a\mapsto\mathfrak m_a$ identifies $K^n$ with the closed points of
$\mathbb A^n_K$, and $K^n\times(K^n)^*$ with the closed points of
$T^*\mathbb A^n_K=\operatorname{Spec}K[x,\xi]$. We write
$X:=K^n$ and $T^*X:=K^n\times(K^n)^*$ for these sets of closed points, with
\emph{affine coordinates} $x=(x_1,\dots,x_n)$. The \emph{projection} is $\pi\colon T^*X\to X$, $(y,\xi)\mapsto y$, and for
$y\in X$ the \emph{fibre} over $y$ is $T_y^*X=\pi^{-1}(y)\cong(K^n)^*$.

The closed subsets of $X$ and of $T^*X$ are the zero sets $V_K(J)$ of ideals
$J\subseteq K[x]$, resp.\ $J\subseteq K[x,\xi]$ (Definition~\ref{def:spec}); by
Theorem~\ref{thm:nss}(a), $V_K(J)$ is the set of closed points of $V(J)$. Via the dual coordinates
$\xi_1,\dots,\xi_n$ we regard polynomials $P\in K[\xi_1,\dots,\xi_n]$ as functions on
$(K^n)^*$; in particular $P(1,0,\dots,0)=P(dx_1)$. 

We write $H:=V_K(x_1)=\{(0,y_2,\dots,y_n)\}\subseteq X$ for the hyperplane on
which the $x_1$-coordinate vanishes; then $\pi^{-1}(H)=H\times (K^n)^* = \{(y,\xi)\in T^*X:y_1 = 0\}\subseteq T^*X$.
The covector $dx_1$, i.e.\ $\xi=(1,0,\dots,0)$, is conormal to $H$: it vanishes
on every direction tangent to $H$.

The projective space $\mathbb P^n$ consists of nonzero vectors
$(q_0,\dots,q_n)\in K^{n+1}$ taken up to a common nonzero scalar (\emph{homogeneous
coordinates}); we embed $X\hookrightarrow\mathbb P^n$ by $x\mapsto(1,x_1,\dots,x_n)$. A
point with $q_0\neq0$ has affine coordinates $x_i=q_i/q_0$, and the points with
$q_0=0$ form the hyperplane at infinity. For a vector $q\in K^{n+1}$ representing a point
of $\mathbb P^n$ we call its entry $i$ its $i$-coordinate, entry $0$ being the
homogenizing coordinate. The \emph{affine part} of a row vector
$\ell=(\ell_0,\dots,\ell_n)$ is $(\ell_1,\dots,\ell_n)$.

\formalref{\leandecl{Stafford38/Paper/Section2/Geometry.lean}{162}{Stafford38.Paper.Section2.proj},
\leandecl{Stafford38/Paper/Section2/Geometry.lean}{166}{Stafford38.Paper.Section2.fibre};}
\begin{definition}[{\cite[VIII.1, V]{lang}}]\label{def:fieldextension}
If $K\subseteq L$ are fields, we call $L$ an extension of $K$ and write $L/K$. For $z\in L$ and $B\subseteq L$ we say:
\begin{enumerate}
\item $z$ is \emph{algebraic over $K$} if there exists a nonzero polynomial $g$
in one variable over $K$ with $g(z)=0$
\item $z$ is \emph{transcendental over $K$} if no such polynomial exists, i.e.\
$g(z)=0\Rightarrow g=0$
\item The \emph{minimal polynomial} of an algebraic $z$ over $K$ is the polynomial of
least degree over $K$ with $g(z)=0$ and coefficient 1 for the highest degree term.
\item $L/K$ is \emph{algebraic} if every element of $L$ is algebraic over $K$
\item $L/K$ is \emph{finitely generated} if $L=K(z_1,\dots,z_m)$ for
  finitely many $z_i$, and \emph{separable} if every $z\in L$ that is algebraic over
  $K$ has a minimal polynomial without repeated roots - in characteristic zero every
  algebraic extension is separable.
\item $B$ is \emph{algebraically independent} over $K$ if
no nonzero polynomial over $K$ in finitely many elements of $B$ vanishes.
\item A maximal algebraically independent subset is a \emph{transcendence basis} of
$L/K$. If $B$ is a transcendence basis, then $L$ is algebraic over the subfield
$K(B)$ generated by $K$ and $B$. 
\item Any two transcendence bases of $L/K$ have the
same cardinality, called the \emph{transcendence degree} $\operatorname{trdeg}_KL$.
\end{enumerate}
The elements of $L$ that are algebraic over $K$ form a subfield, and every
finite set generating $L$ over $K$ contains a transcendence basis.
\end{definition}

\begin{definition}[{\cite[§25]{matsumura}}]\label{def:omega}
Let $\kappa/K$ be a field extension.
\begin{enumerate}[label=(\alph*)]
\item A \emph{$K$-derivation} of $\kappa$ is a $K$-linear map
$D\colon\kappa\to\kappa$ with $D(ab)=aD(b)+bD(a)$.
\item Let $\mathcal G$ be the free $\kappa$-vector space with basis symbols
$\langle a\rangle$, $a\in\kappa$. The space of \emph{differentials}
$\Omega_{\kappa/K}$ is the quotient of $\mathcal G$ by the subspace spanned by
\[
 \langle a+b\rangle-\langle a\rangle-\langle b\rangle,\quad
 \langle ab\rangle-a\langle b\rangle-b\langle a\rangle,\quad
 \langle c\rangle\qquad(a,b\in\kappa,\ c\in K),
\]
and $da\in\Omega_{\kappa/K}$ denotes the class of $\langle a\rangle$.
\item A $K$-derivation $D$ induces the $\kappa$-linear map
$\Omega_{\kappa/K}\to\kappa$, $da\mapsto D(a)$, since the linear map
$\langle a\rangle\mapsto D(a)$ on $\mathcal G$ kills the relations. Conversely, a
$\kappa$-linear map $\Omega_{\kappa/K}\to\kappa$ gives the $K$-derivation
$a\mapsto(\text{its value on }da)$. Hence a basis $d\beta_1,\dots,d\beta_d$ of
$\Omega_{\kappa/K}$ has dual derivations $D_1,\dots,D_d$ with
$D_k(\beta_l)=\delta_{kl}$.
\end{enumerate}
\end{definition}
\begin{definition}[{\cite[§10-§11]{matsumura}}]\label{def:valuation}
Let $F$ be a field. 
\begin{enumerate}[label=(\alph*)]
\item A subring $W\subseteq F$ is a \emph{valuation ring} of $F$ if for
every $a\in F^\times$ we have $a\in W$ or $a^{-1}\in W$. It is a local ring with
maximal ideal $\mathfrak m_W=W\setminus W^\times$. 
\item A \emph{discrete valuation ring} (DVR) is a Noetherian valuation ring $V\ne F$. Then $\mathfrak m_V$ is generated by some $\varpi\in V$, and every $f\in F^\times$ is $f=\varpi^ku$ with $k\in\mathbb Z$
and $u\in V^\times$, and $\operatorname{ord}_V(f):=k$. 
\item For a field $\kappa = V/\mathfrak m_V$ we write $\kappa[[t]]$ for the ring of formal power series, and $\kappa((t))$ for its fraction field (formal Laurent series).
\end{enumerate}
\end{definition}
\begin{definition}\label{def:conormal}
We say
\begin{enumerate}[label=(\alph*)]
\item A closed set $Y$ is \emph{irreducible} if it is not a union of two proper closed subsets, equivalently if $I(Y)$ is prime \cite[Cor I.1.4]{hartshorne}. 
\item For irreducible closed $Y$, the \emph{smooth locus} $Y_{\mathrm{sm}}$ is the set of $y\in Y$ at which the Jacobian matrix of generators of $I(Y)$ has rank $n-\dim Y$, i.e.\ the points where $Y$ has a well-defined tangent space of the expected dimension $\dim Y$ (i.e. without kinks etc.). It is open and dense in $Y$. \cite[Ch. I.5]{hartshorne} 
\item For $y\in Y$ the \emph{conormal space} is
\[N^*_yY=\operatorname{span}_K\{\,df_y : f\in I(Y)\,\}\subseteq(K^n)^*,\qquad
 df_y=\textstyle\sum_i\partial_if(y)\,dx_i,\]
 \item $N(Y)=\bigcup_{y\in Y_{\mathrm{sm}}}N^*_yY\subseteq (K^n)^*$ denotes the \emph{conormal
cone}. 
\end{enumerate}
\end{definition}

\begin{definition}\label{def:dir}
For $Y$ as above, the set of \emph{conormal directions} is
\[ \operatorname{Dir}(Y)=\overline{\{[\xi]\in\mathbb P((K^n)^*): 0\neq\xi\in N(Y)\}}
 \subseteq\mathbb P((K^n)^*),\]
the Zariski closure of the projectivization of $N(Y)$.
\end{definition}

\begin{definition}\label{def:fibre-can}
A closed subset $W\subseteq T^*X$ is \emph{fibre-conical} if it is stable under the
fibrewise scaling $(y,\xi)\mapsto(y,t\xi)$ for all $t\in K^\times$.\\
We say it is \emph{involutive} if $I(W)$ is.
\end{definition}

\begin{theorem}\label{thm:asymptotic-conormal-formal}
Let $Y\subseteq X$ be nonempty, closed, and irreducible with $Y\cap H=\varnothing$, and let
$N(Y)=\bigcup_{y\in Y_{\mathrm{sm}}}N^*_yY\subseteq (K^n)^*$ be the conormal cone over the smooth locus $Y_{sm}$.
\begin{enumerate}[label=(\roman*)]
\item Every polynomial $P\in K[\xi_1,\dots,\xi_n]$ vanishing on $N(Y)$ satisfies $P(1,0,\dots,0)=0$.
\item $[dx_1]\in\operatorname{Dir}(Y)$.
\end{enumerate}
\end{theorem}

\formalref{%
\leandecl{Stafford38/Paper/Section5/SmoothLocus.lean}{645}{Stafford38.Paper.Section5.theorem_5_7_i_printed},
\leandecl{Stafford38/Paper/Section5/SmoothLocus.lean}{660}{Stafford38.Paper.Section5.theorem_5_7_ii_printed}}

\begin{theorem}\label{thm:coisotropic-exclusion}
Let $K$ be an algebraically closed field of characteristic zero, and let
$W\subseteq T^*X$ be a Zariski closed subset that is fibre-conical and
involutive. Let $P\in K[\xi_1,\dots,\xi_n]$ be
a polynomial with $P(1,0,\dots,0)\neq0$. If $P$ vanishes on $W$ and $W\cap\pi^{-1}(H)=\varnothing$, then $W$ must be empty.
\end{theorem}

\formalref{%
\leandecl{Stafford38/Paper/Section5/Exclusion.lean}{447}{Stafford38.Paper.Section5.theorem_5_8}}

\subsection{Proof of Theorem \ref{thm:asymptotic-conormal-formal}}\label{sub:asym}
\begin{proof}

If $Y$ is nonempty, closed and irreducible, then $\mathfrak p:=I(Y)$ is a prime
ideal and $Y=V_K(\mathfrak p)$. We will write $R=K[x]/\mathfrak p$, $F=\operatorname{Frac}(R)$ for the field of fractions, and $\bar x_i$ for the image of $x_i$ in $R$. \\

If $\bar x_1$ is algebraic over $K$, there exists a nonzero polynomial $g$ over $K$
in one variable with $g(\bar x_1)=0$ (Def. \ref{def:fieldextension}), therefore \ $g(x_1)\in\mathfrak p$. Since $K$
is algebraically closed, $g(x_1)=\lambda\prod_i(x_1-c_i)$ with
$\lambda\in K^\times$, and as $\mathfrak p$ is prime, some factor $x_1-c$ lies
in $\mathfrak p$. Thus $x_1=c$ on $Y=V_K(\mathfrak p)$, and since $Y\cap H = \varnothing$ and $Y\neq\varnothing$, we have $c\neq0$.
In this case $f=x_1-c$ lies in $\mathfrak p=I(Y)$ and has the differential
$df_y=dx_1$ at every point $y\in Y_{\mathrm{sm}}$, so $dx_1\in N_y^*Y$, and its coordinate vector in the dual basis $(1,0,\dots , 0) \in N(Y)$. Every polynomial vanishing on $N(Y)$ therefore vanishes at $(1,0,\dots , 0)$, and
$[dx_1]\in\operatorname{Dir}(Y)$. So we may assume from now on that
$\bar x_1$ is transcendental over $K$.\\

The rest of the proof follows this rough logic. We construct a discrete valuation
ring $V$ of $F$ at infinity in which $x_1$ vanishes, and the corresponding
generic point $y\in\kappa((t))^n$ of $Y$. The affine cone $\mathcal C$ of $Y$ over the projective tangent plane at
$y$ has a basis $M$ of power series in $t$ that remains a basis at $t=0$, and
whose reduction $M(0)$ has vanishing $1$-coordinate. Hence the row vector
$u=(0,1,0,\dots,0)$, which reads off the $1$-coordinate, satisfies
$uM(0)=0$, and we lift it to row vectors $\ell(t)$ annihilating $\mathcal C$,
with $\ell(0)=u$. By construction the affine part $\tilde\ell(t)$ consists of
conormal covectors to $Y$ at $y$, and $\tilde\ell(0)=dx_1$. A polynomial
vanishing on $N(Y)$ therefore vanishes at $dx_1$, and so $[dx_1]\in\operatorname{Dir}(Y)$.

We set $\bar x_0:=1$.

\begin{lemma}\label{lem:W}
Let $\bar x_1$ be transcendental over $K$. There exist a valuation ring $W$ with
$K\subseteq W\subset F$, $\bar x_1\in\mathfrak m_W$, and an
index $j$ with $\bar x_i/\bar x_j\in W$ for all $0\le i\le n$.
\end{lemma}
\begin{proof}
Since $\bar x_1$ is transcendental over $K$, $S=K[\bar x_1]$ is a polynomial
ring, and $\mathfrak q=\bar x_1S$ is a prime ideal of $S$. The localization
$S_{\mathfrak q}\subseteq F$ is a local ring with maximal ideal
$\bar x_1S_{\mathfrak q}$. By \cite[Tag~00IA]{Stacks} there exists a valuation ring
$W$ of $F$ dominating $S_{\mathfrak q}$, i.e.\ $S_{\mathfrak q}\subseteq W$ and
$\mathfrak m_{S_{\mathfrak q}}=\mathfrak m_W\cap S_{\mathfrak q}$. Hence $K\subseteq W$ and
$\bar x_1\in\mathfrak m_{S_{\mathfrak q}}\subseteq\mathfrak m_W$. For nonzero $a,b\in F$ either $a/b\in W$ or $b/a\in W$,
so among the finitely many nonzero $\bar x_i$ we can choose $j$ with
$\bar x_i/\bar x_j\in W$ for all $i$ (note $\bar x_0=1\ne0$).
\end{proof}

Let $W$ be the valuation ring of Lemma \ref{lem:W}. We define the column vector $q$ in homogeneous coordinates by
\[ q_i:=\bar x_i/\bar x_j\qquad(0\le i\le n),\]
so that $q_i\in W$, $q_j=1$, $q_0=\bar x_j^{-1}\neq0$, $\bar x_i=q_i/q_0$ and
$q_1=q_0\bar x_1$.

\begin{lemma}\label{lem:place}
For $q$ defined as above, there exists a discrete valuation ring $V$ with $K\subseteq V\subset F$,
$\operatorname{Frac}V=F$, such that, with $\kappa=V/\mathfrak m_V$ and $\operatorname{res}(q_i)$ the
residue of $q_i$ in $\kappa$:
\begin{enumerate}[label=(\alph*)]
\item $q_i\in V$ for all $i$, $\bar x_1\in\mathfrak m_V$ and $R\not\subseteq V$;
\item $\kappa$ is algebraic over $L=K(\operatorname{res}(q_0),\dots, \operatorname{res}(q_n))$;
\item $\Omega_{V/K}$ is a finitely generated $V$-module.
\end{enumerate}
\end{lemma}
\begin{proof}
Let $C_0=K[q_0,\dots,q_n,\bar x_1]\subseteq W$; it is finitely generated with
$\operatorname{Frac}C_0=F$ since on one hand $\bar x_i = q_i/q_0$ so every $\bar x_i$ is in $\operatorname{Frac}C_0$ and on the other hand, $F$ is generated by the $\bar x_i$. Let $A$ be its integral closure in $F$. Since $C_0$ is a domain finitely generated over a field, $A$ is a finite $C_0$-module \cite[Thm. I.3.9A]{hartshorne}, hence a finitely generated, integrally
closed Noetherian domain with $\operatorname{Frac}A=F$ and $q_i\in A$. As $W$ is a valuation ring and thereby is
integrally closed \cite[Prop 5.18]{atiyah}, $A\subseteq W$. $\bar x_1 \in \mathfrak m_W$, so $x_1$ is a nonzero nonunit of $W$ and therefore also of $A$. Since $\bar x_1$ is a nonzero nonunit of the Noetherian domain $A$, the ideal
$\bar x_1A$ is proper, so some prime $\mathfrak p_D$ is minimal among the primes
containing $\bar x_1$. By Krull's principal ideal theorem
\cite[Corr. 11.17]{atiyah}, $\operatorname{ht}\mathfrak p_D\le1$. Because $A$
is a domain, the only prime of height $0$ is $(0)$, and $\bar x_1\neq0$ lies in
$\mathfrak p_D$; hence $\operatorname{ht}\mathfrak p_D=1$.

Put $V:=A_{\mathfrak p_D}$. This is a Noetherian local domain. It is normal,
since a localization of an integrally closed domain is integrally closed \cite[Prop 5.12]{atiyah}, and
its dimension equals $\operatorname{ht}\mathfrak p_D=1$. A Noetherian local
normal domain of dimension one is a discrete valuation ring
\cite[Thm. 11.2]{matsumura}, so $V$ is a discrete valuation ring.

It remains to check the listed properties.
\begin{itemize}
\item $K\subseteq A\subseteq V$: we have $K\subseteq C_0\subseteq A$, and the
domain $A$ embeds into its localization $V$.
\item $\operatorname{Frac}V=F$: since $A\subseteq V\subseteq F$ and
$\operatorname{Frac}A=F$.
\item $q_i\in V$ for all $i$: since $q_i\in A\subseteq V$.
\item $\bar x_1\in\mathfrak m_V$: since $\mathfrak m_V=\mathfrak p_DA_{\mathfrak p_D}$
and $\bar x_1\in\mathfrak p_D$. $V_K(p, x_1) = Y \cap H = \varnothing$, so $1 \in (p, x_1)$ by Theorem \ref{thm:nss}, i.e. $1=p+gx_1$ for some $p\in\mathfrak p$, $g\in K[x]$. Modulo
$\mathfrak p$ this gives $1=\bar g\,\bar x_1$, so $\bar x_1\in R^\times$.
If $R\subseteq V$, then $\bar x_1\in R^\times\subseteq V^\times$, contradicting
$\bar x_1\in\mathfrak m_V$.
\end{itemize}
This proves (a).

(b) Since $V=A_{\mathfrak p_D}$ and $\mathfrak m_V=\mathfrak p_DV$, we have
$\kappa=V/\mathfrak m_V=(A/\mathfrak p_D)_{(0)}=\operatorname{Frac}(A/\mathfrak p_D)$,
and the residue map $\operatorname{res}:V\to\kappa$ restricts to the quotient map $A\to A/\mathfrak p_D$. As $A$ is integral over $C_0$ and $\bar x_1\in\mathfrak p_D$, the ring $A/\mathfrak p_D$ is integral over the image of
$C_0$ in $A_{\mathfrak p_D}$, which is
$K[\operatorname{res}(q_0),\dots,\operatorname{res}(q_n)]\subseteq L$. So every element of
$A/\mathfrak p_D$, and hence of its fraction field $\kappa$, is algebraic over $L$.

(c) $A$ is a finitely generated $K$-algebra, so $\Omega_{A/K}$ is a finite
$A$-module, and $\Omega_{V/K}=V\otimes_A\Omega_{A/K}$ since $V$ is a localization of $A$ \cite[Prop II.8.2A]{hartshorne}.
\end{proof}

By Lemma~\ref{lem:place} there exists a discrete valuation ring $V\subset F$ with
$q_i\in V$ for all $i$, $\bar x_1\in\mathfrak m_V$ and $R\not\subseteq V$. Let $\varpi$ be a
generator of $\mathfrak m_V$, then every nonzero $f \in F$ can be written as $f = \varpi^k\cdot u$ with $k \in \mathbb{Z}$ and $u$ a unit of $V$ - we say $\operatorname{ord}_V(f) = k \in \mathbb Z$ is the valuation of $f$ in $V$. We denote the residual field with $\kappa=V/\mathfrak m_V$.

Let $\widehat V$ be the completion of $V$. Since $V\supseteq K$
and $K$ is of characteristic zero, there exists a coefficient field of $\widehat V$
containing $K$ and mapping isomorphically onto $\kappa$ \cite[Thm. 28.3]{matsumura};
choosing one gives us an isomorphism $\widehat V\cong\kappa[[t]]$ under which
$\varpi$ corresponds to the formal variable $t$. Since the natural map $V\to\widehat V$ is injective ($\bigcap_k\varpi^kV=0$ \cite[Corr. 10.20]{atiyah}) and extends to an embedding $F\hookrightarrow\kappa((t))$, we can identify elements of $F$ with
Laurent series. The residue $\operatorname{res}(f)$ of $f \in V$ is then just the constant term $f(0)$ of its series; in particular $\operatorname{res}(q_i) = q_i(0).$

Let $a=\operatorname{ord}_V(q_0)$ and $b=\operatorname{ord}_V(q_1)$. $a,b\geq 0$ since $q_0,q_1$ are elements of $V$, so that in $\kappa[[t]]$
\[q_0(t)=t^au_0(t),\qquad q_1(t)=t^bu_1(t),\]
with $u_0(0)\neq0\neq u_1(0)$. Since $\bar x_1\in\mathfrak m_V$, $b-a=\operatorname{ord}_V(q_1/q_0) = \operatorname{ord}_V(\bar x_1)>0$. Moreover
$a>0$: if $a=0$, then $q_0$ is a unit of $V$, hence all $\bar x_i=q_i/q_0$
lie in $V$, so $R\subseteq V$ since also $K\subseteq V$, contradicting Lemma~\ref{lem:place}(a). It follows that $j \notin \{0,1\}$ since $q_j = 1$ is a unit while $q_0,q_1$ are not.

Dividing by $q_0(t)\neq 0$, we define $y=(q_1(t)/q_0(t),\dots,q_n(t)/q_0(t))\in\kappa((t))^n$ in affine coordinates. Its
coordinates are the images of $\bar x_1,\dots,\bar x_n$ under $F\hookrightarrow
\kappa((t))$. Since ring homomorphisms commute with polynomial evaluation and
this map is injective, for every polynomial $f\in K[x_1,\dots,x_n]$ we have
\begin{equation}
    f(y)=0\iff f(\bar x_1,\dots,\bar x_n)=0 \iff  f\in\mathfrak p.
    \label{eq:generic}
\end{equation}

Let $df_y=\nabla f(y)\in\kappa((t))^n$ for a polynomial $f$.
We define $\Theta_y$ to be the tangent space of $Y$ at the generic point $y$,
\[ \Theta_y=\bigcap_{f\in\mathfrak p}\ker df_y\subseteq\kappa((t))^n,\]
the directions $v$ with $df_y\cdot v=0$ for all equations $f$ of $Y$. We now pass to homogeneous coordinates by defining a map
\[ \rho\colon\kappa((t))^{n+1}\to\kappa((t))^n,\qquad
 \rho(v(t))_i=\frac{v_i(t)q_0(t)-q_i(t)v_0(t)}{q_0^2(t)}\quad(1\le i\le n),\]
which sends a homogeneous vector $v(t)$ to the affine direction it induces at $q(t)$
(by the quotient rule); its kernel is the line $\kappa((t))q(t)$. We write
$\mathcal C=\rho^{-1}(\Theta_y)$ for the affine cone over the projective tangent
plane. Our goal is to build a basis of $\mathcal C$ whose entries are power
series in $t$ and which also remains a basis at $t=0$. 
We do so by building independent vectors from the tangent planes of $q(t)$.
Since $0<a<b$, in $\kappa$ we have $q_0(0)=q_1(0)=0$ and $q_j(0)=1$. Let
$L=K(q_0(0),\dots,q_n(0))\subseteq\kappa$; by Lemma~\ref{lem:place}(b),
$\kappa/L$ is algebraic.

\begin{lemma}\label{lem:kappa}
There exist indices $i_1,\dots,i_{d'}\notin\{0,1,j\}$ such that
$d(q_{i_1}(0)),\dots,d(q_{i_{d'}}(0))$ form a basis of $\Omega_{\kappa/K}$. In
particular $d' = \dim_\kappa\Omega_{\kappa/K}<\infty$.
\end{lemma}
\begin{proof}

Since $q_0(0)=q_1(0)=0$ and $q_j(0)=1$, the field $L$ is generated over $K$ by
the finitely many $q_i(0)$ with $i\notin\{0,1,j\}$, so some subset
$B=\{q_{i_1}(0),\dots,q_{i_{d'}}(0)\}$ of them is a transcendence basis of $L/K$
\cite[Thm. VIII.1.1]{lang}. As $L/K(B)$ and $\kappa/L$ are algebraic
(Lemma~\ref{lem:place}(b)), so is $\kappa/K(B)$; hence $B$ is also a
transcendence basis of $\kappa/K$. Since $K$ has characteristic zero, a
transcendence basis of $\kappa/K$ is a differential basis, i.e.\ its
differentials form a $\kappa$-basis of $\Omega_{\kappa/K}$
\cite[Thm.~26.5]{matsumura}.
\end{proof}

We write for all $q_i(t) = \sum_rc_{ir}t^r$. The derivative with respect to $t$ is then $q'_i(t) = \sum_rr\,c_rt^{r-1}$. This is a $K$-derivation of $\kappa((t))$. By Lemma \ref{lem:kappa} there exist $i_1,\dots,i_{d'}\notin\{0,1,j\}$ such that $\{d(q_{i_1}(0)),\dots,d(q_{i_{d'}}(0))\}$ form a basis of $\Omega_{\kappa/K}$. Let $D_1,\dots,D_{d'}$ be the dual $K$-derivations of this basis (Definition~\ref{def:omega}(c)). 
We apply $D_k$ coefficientwise to $q_i(t)=\sum_rc_{ir}t^r$ and define $Z_{ik}(t):=\sum_rD_k(c_{ir})t^r$. Since $D_k(q_{i_l}(0))=\delta_{kl}$, the reduction $Z(0)$ is the identity on the rows $i_1,\dots,i_{d'}$, and since $q_j(t) = 1$ is constant, $Z_{jk}(t)=0$. Rows $0$ and $1$ of $Z(t)$ have order at least $a$ and $b$, since $D_k$ does not differentiate with respect to t. 

Let $\tilde Z(t)$ be the ${d'}\times {d'}$ submatrix of $Z(t)$ consisting of the rows
$i_1,\dots,i_{d'}$ and all columns. Its reduction $\tilde Z(0)$ is the identity
matrix, so $\tilde Z(t)$ is invertible over $\kappa[[t]]$. Hence there exists a unique vector
$\lambda(t)=(\lambda_1(t),\dots,\lambda_{d'}(t))^{\mathsf T}\in\kappa[[t]]^{d'}$ with
$\tilde Z(t)\lambda=(q'_{i_1}(t),\dots,q'_{i_{d'}}(t))^{\mathsf T}$. We define $w(t)=q'(t)-Z(t)\lambda(t)$, so
that $w_{i_l}(t)=0$ for $l=1,\dots,{d'}$.
Its $0$th entry is $w_0(t)=q_0'(t)-Z_{0\bullet}(t)\lambda(t)$. Here
$q_0'(t)=at^{a-1}u_0(t)+t^au_0(t)'$ has order exactly $a-1$ (since the coefficient
$au_0(0)$ of $t^{a-1}$ is nonzero in characteristic zero),
while $Z_{0\bullet}(t)\lambda(t)$ has order at least $a$, because every entry of
row $0$ of $Z(t)$ has order at least $a$ (since $q_0(t)$
has order $a$). Hence $\operatorname{ord}w_0(t)=a-1$. Let $c$ be the minimal order
among the entries of $w(t)$; then $c\le a-1$, and $\tau(t):=t^{-c}w(t)$ is a vector
over $\kappa[[t]]$ with $\tau(0)\neq0$ and $\tau_{i_l}(t)=0$. In row $1$, both
$q_1'(t)$ and $Z_{1\bullet}(t)\lambda(t)$ have order at least $b-1$, so $\operatorname{ord}w_1(t)\ge b-1$, so
$\operatorname{ord}\tau_1(t)\ge b-1-c\ge b-a>0$. Finally $q_j'(t)=0$ and
$Z_{j\bullet}(t)=0$ (as $q_j(t)=1$), so $\tau_j(t)=0$. \\

We now show that $q(t), Z_{\bullet1}(t),\dots, Z_{\bullet {d'}}(t), \tau(t)$ span all of $\mathcal C$: 
Let $M(t)=(\,q(t), Z_{\bullet1}(t),\dots, Z_{\bullet {d'}}(t),\tau(t)\,)$, an
$(n+1)\times({d'}+2)$ matrix over $\kappa[[t]]$. Its reduction $M(0)$ has full
column rank: in row $j$ only $q(0)$ is nonzero ($q_j(0)=1$); in rows $i_1,\dots,i_{d'}$
the columns $Z_{\bullet k}(0)$ form an identity block while $\tau(0)$ vanishes
there; and $\tau(0)\neq0$, i.e. we have ${d'}+2$ independent reduced column vectors. This also translates to the column vectors of $M(t)$ themselves: Since some $({d'}+2)\times({d'}+2)$ minor of $M(0)$ is nonzero, the corresponding
minor of $M(t)$ is a unit of $\kappa[[t]]$; hence the columns of $M(t)$ are also
linearly independent over $\kappa((t))$, so they span a $({d'}+2)$-dimensional
subspace there.
Furthermore, every column of $M(0)$ has a vanishing entry for the $1$-coordinate:
$q_1(0)=0$ (as $b>0$), the entries $Z_{1k}(t)$ have order at least $b>0$, hence $Z_{1k}(0) = 0$, and
$\tau_1(0)=0$. Thus the span of $q(0), Z_{\bullet1}(0),\dots, Z_{\bullet {d'}}(0), \tau(0)$ is contained in $ \{v \in \kappa^{n+1} : v_1 = 0\}$ and also has dimension ${d'}+2$. 

To show that the columns of $M(t)$ really span all of $\mathcal C$, we must still show that a) the dimensions match and b) the columns actually lie in $\mathcal C$. The following Lemma will help us with that.
\begin{lemma}\label{lem:tangent}
Let $m = \operatorname{trdeg_K}F$ and  $\Theta_y$, $\rho$ and $\mathcal C$ be defined as above.
\begin{enumerate}[label=(\alph*)]
\item $\dim_{\kappa((t))}\Theta_y = m \le {d'}+1$; hence $\dim_{\kappa((t))}\mathcal C = m + 1 \le {d'}+2$.
\item For every $K$-derivation $\delta$ of $\kappa((t))$ the column
$\delta q(t)=(\delta q_i(t))_i$ lies in $\mathcal C$; moreover $q(t)\in\mathcal C$.
\end{enumerate}
\end{lemma}
\begin{proof}
(a) For a finitely generated field extension $F/K$ in characteristic $0$, $\dim_F\Omega_{F/K}=\operatorname{trdeg}_KF = m$ \cite[Thm.~II.8.6A]{hartshorne}.
The second fundamental exact sequence \cite[Theorem 25.2]{matsumura} for $K[x]\to R=K[x]/\mathfrak p$ is
\[\mathfrak p/\mathfrak p^2\to\Omega_{K[x]/K}\otimes_{K[x]}R\to\Omega_{R/K}\to0 ,\]
where $\Omega_{K[x]/K}$ is free on $dx_1,\dots,dx_n$, so the middle term is $R^n$
and the first map sends $f$ to $df(\bar x)$. Hence
$\Omega_{R/K}=R^n/\operatorname{span}_R\{df(\bar x):f\in\mathfrak p\}$.
For a multiplicative set $S\subseteq R$, $\Omega_{S^{-1}R/K}=S^{-1}\Omega_{R/K}$ \cite[Prop II.8.2A]{hartshorne}. In particular, this means for $F=S^{-1}R$ with
$S=R\smallsetminus\{0\}$ $\Omega_{F/K}=F\otimes_R\Omega_{R/K}$ and since localization is exact, applying $F\otimes_R-$ to the description of
$\Omega_{R/K}$ above gives
\[ \Omega_{F/K}=F^n/\operatorname{span}_F\{df(\bar x):f\in\mathfrak p\}.\]
Analogously, applying this for $S = V\setminus \{0\}$ and $F=S^{-1}V$ gives $F\otimes_V\Omega_{V/K}=\Omega_{F/K}$.
$\operatorname{span}_F\{df(\bar x):f\in\mathfrak p\}$ has therefore dimension $n-m$,
i.e.\ the Jacobian $(\partial_if(\bar x))_{f\in\mathfrak p}$ has rank $n-m$ over $F$.
Rank does not change under the field extension $F\subseteq\kappa((t))$, so
$\dim\operatorname{span}_{\kappa((t))}\{df_y:f\in\mathfrak p\}=n-m$, and therefore
$\dim_{\kappa((t))}\Theta_y=\dim_{\kappa((t))}\bigcap_f\ker df_y=n-(n-m)=m$.

The second fundamental exact sequence for $V\to\kappa$ is
\[ \mathfrak m_V/\mathfrak m_V^2\to\kappa\otimes_V\Omega_{V/K}\to\Omega_{\kappa/K}\to0 .\]
Since $\dim_\kappa\mathfrak m_V/\mathfrak m_V^2=1$ and
$\dim_\kappa\Omega_{\kappa/K}=d'$, we get
$\dim_\kappa(\kappa\otimes_V\Omega_{V/K})\le d'+1$. As $\Omega_{V/K}$ is a finitely
generated $V$-module (Lemma~\ref{lem:place}(c)), Nakayama's lemma shows that it is
generated by at most $d'+1$ elements. Since $F\otimes_V\Omega_{V/K}=\Omega_{F/K}$,
we obtain $m\le d'+1$, hence $\dim\Theta_y\le d'+1$. Finally, $\rho$ maps
$\mathcal C$ onto $\Theta_y$ with one-dimensional kernel, so
$\dim_{\kappa((t))}\mathcal C=\dim_{\kappa((t))}\Theta_y+1\le d'+2$.

(b) Let $f\in\mathfrak p$ and $f = \sum_\alpha c_\alpha y^\alpha$. By \eqref{eq:generic}, $f(y)=0$, and $\delta$ vanishes on
$K$, so
\[\delta\bigl(f(y)\bigr)
=\sum_\alpha\Bigl(\underbrace{y^\alpha\,\delta(c_\alpha)}_{=\,0}+c_\alpha\,\delta(y^\alpha)\Bigr)
=\sum_{i=1}^n\frac{\partial f}{\partial x_i}(y)\,\delta y_i
=df_y\cdot\delta y,\] hence
$\delta y\in\Theta_y$. Since $y_i=q_i/q_0$, the quotient rule gives
$\rho(\delta q)_i=(\delta q_i\,q_0-q_i\,\delta q_0)/q_0^2=\delta(q_i/q_0)=\delta y_i$.
Thus $\rho(\delta q)=\delta y\in\Theta_y$, i.e.\ $\delta q\in\mathcal C$. Moreover, $\rho(q)_i=(q_iq_0-q_iq_0)/q_0^2=0$ for all $i$, so $\rho(q)=0\in\Theta_y$
and therefore $q\in\mathcal C$.

\end{proof}

The column vectors $q(t)$, $Z_{\bullet k}(t)$ and $q'(t)$ are of the form $\delta q(t)$, where
$\delta$ is a $K$-derivation of $\kappa((t))$ (the $D_k$ extended
coefficientwise, and $d/dt$). By Lemma~\ref{lem:tangent}(b) they lie in $\mathcal C$. Since
$\tau(t)=t^{-c}(q'(t)-Z(t)\lambda(t))$ and $t^{-c}\in\kappa((t))$, $\tau(t)$ is a
$\kappa((t))$-linear combination, hence also lies in $\mathcal C$. By Lemma~\ref{lem:tangent}(a), $\dim \mathcal C\le {d'}+2$, so the ${d'}+2$ independent columns of $M(t)$ must span $\mathcal C$.\\

\begin{lemma}\label{lem:ell}
If a row vector $\ell\in\kappa((t))^{n+1}$ annihilates $\mathcal C$, then
$(\ell_1,\dots,\ell_n)\in\operatorname{span}\{df_y:f\in\mathfrak p\}$.
\end{lemma}
\begin{proof}
Let $\ell$ annihilate $\mathcal C$. Since $q\in\mathcal C$, $\ell(q)=0$, i.e.\
$\sum_{i\ge1}\ell_iq_i=-\ell_0q_0$. For $v\in\mathcal C$ we have $\ell(v)=0$, i.e.\
$\sum_{i\ge1}\ell_iv_i=-\ell_0v_0$, so
\[ \sum_{i\ge1}\ell_i\rho(v)_i
 =\frac{q_0\sum_{i\ge1}\ell_iv_i-v_0\sum_{i\ge1}\ell_iq_i}{q_0^2}
 =\frac{-q_0\ell_0v_0+v_0\ell_0q_0}{q_0^2}=0 .\]
As $\rho(\mathcal C)=\Theta_y$, the vector $(\ell_1,\dots,\ell_n)$ annihilates
$\Theta_y=\bigcap_f\ker df_y$, hence lies in
$\Theta_y^\perp=\operatorname{span}\{df_y:f\in\mathfrak p\}$.
\end{proof}

Let $u=(0,1,0,\dots,0)$ be the row vector reading off the $1$-coordinate;
since $\operatorname{im}M(0)\subseteq \{v \in \kappa^{n+1} : v_1 = 0\}$, we know that $uM(0)=0$.
We can lift $u$ to a row vector $\ell(t)\in\kappa[[t]]^{n+1}$ with $\ell(t) M(t)=0$ and
$\ell(0)=u$: Since $M(0)$ has full column rank ${d'}+2$, there is a set $\iota$ of
${d'}+2$ row indices for which the square submatrix $M_\iota(0)$ of $M(0)$ is
invertible; then $M_\iota$ is invertible over $\kappa[[t]]$, since its determinant
has a nonzero constant term. Let $\upsilon$ be the complementary set of row indices,
and define $\ell(t)$ by
\[\ell_\upsilon(t)=u_\upsilon,\qquad \ell_\iota(t)=-u_\upsilon M_\upsilon (t)M_\iota^{-1}(t),\]
where $M_\upsilon(t)$ denotes the submatrix of $M(t)$ with rows in $\upsilon$. Then
$\ell(t) M(t)=\ell_\iota(t)M_\iota(t)+\ell_\upsilon(t)M_\upsilon(t)=-u_\upsilon M_\upsilon(t)+u_\upsilon M_\upsilon(t)=0$. Moreover
$uM(0)=0$ gives $u_\iota=-u_\upsilon M_\upsilon(0)M_\iota(0)^{-1}$, hence $\ell(0)=u$.

Since the columns of $M(t)$ span $\mathcal C$, $\ell(t)$ annihilates $\mathcal C$, so by
Lemma~\ref{lem:ell} its affine part
$\tilde\ell(t)=(\ell_1(t),\dots,\ell_n(t))$ lies in
$\operatorname{span}\{df_y:f\in\mathfrak p\}$, and since $\ell(0)=u$,
$\tilde\ell(0)=(1,0,\dots,0)$.

Finally, let $P\in K[\xi_1,\dots,\xi_n]$ vanish on $N(Y)$. We will show
$P(\tilde\ell(0))=0$: Since $\tilde\ell(t)\in\operatorname{span}_{\kappa((t))}\{dg_y:g\in
\mathfrak p\}$, we can write
\[\tilde\ell(t)=\sum_{k=1}^s\mu_k(t)\nabla g_k(y),\qquad g_k\in\mathfrak p,\
 \mu_k(t)\in\kappa((t)),\]
 With some $s \in \mathbb{N}$. Here $\nabla g(x)=(\partial g/\partial x_1(x),\dots,
\partial g/\partial x_n(x))$ denotes the coordinate vector of the covector
$dg_x$. We furthermore define
\[\Phi(x,\mu)=P\Bigl(\sum_{k=1}^s\mu_k\nabla g_k(x)\Bigr)
 =\sum_\alpha\phi_\alpha(x)\mu^\alpha\in K[x_1,\dots,x_n,\mu_1,\dots,\mu_s],\]
with $\phi_\alpha\in K[x_1,\dots,x_n]$ and $\mu^\alpha = \mu_1^{\alpha_1}\dots\mu_s^{\alpha_s}$ for a multi-index $\alpha \in \mathbb{N}^s$. 
For a smooth point $y'\in Y_{\mathrm{sm}}$ and
$\mu\in K^s$ we have $\sum_k\mu_kdg_k(y')\in N^*_{y'}Y\subseteq N(Y)$, so
$\Phi(y',\mu)=0$. As $K$ is infinite, a polynomial in $\mu$ vanishing on all of
$K^s$ is zero, so $\phi_\alpha(y')=0$ for all $\alpha$. Since $Y_{\mathrm{sm}}$
is dense in $Y$, every $\phi_\alpha$ vanishes on all of $Y$, i.e.\
$\phi_\alpha\in\mathfrak p=I(Y)$. By \eqref{eq:generic} (applied to the generic
point $y\in\kappa((t))^n$), we therefore have $\phi_\alpha(y)=0$ for all
$\alpha$. Hence $\Phi(y,\mu(t))=\sum_{\alpha}\phi_\alpha(y)\mu^\alpha(t)=\sum_\alpha 0\cdot \mu^\alpha(t) = 0$ for all $\mu(t)\in\kappa((t))^s$, and in
particular
\[ P(\tilde\ell(t))=\Phi\bigl(y,(\mu_1(t),\dots,\mu_s(t))\bigr)=0.\]
Now $\tilde\ell(t)\in\kappa[[t]]^n$ and $P$ has coefficients in $K$, so
$P(\tilde\ell(t))\in\kappa[[t]]$, and its constant term is $P(\tilde\ell(0))$. Since $P(\tilde\ell(t))=0$, we get
$P(\tilde\ell(0))=P(1,0,\dots,0)=0$, which is statement~(i).

By statement~(i), every polynomial vanishing on $N(Y)$ vanishes at $(1,0,\dots,0)$. We
use this to show $[dx_1]\in\operatorname{Dir}(Y)$. Recall that $\operatorname{Dir}(Y)$ is the Zariski closure of the set of classes $[\xi]$ with $0\ne\xi\in N(Y)$, and that a point of $\mathbb P((K^n)^*)$ lies in the closure of a set $S$ if and only if every homogeneous polynomial vanishing on $S$ also vanishes at that point. So let $G$ be a homogeneous polynomial vanishing on all classes $[\xi]$ with $0\ne\xi\in N(Y)$. Since also $G(t\xi)=t^{\deg G}G(\xi)$ by homogeneity, $G$ vanishes on every nonzero $\xi\in N(Y)$ (and at $0$ if $\deg G>0$), hence on $N(Y)$; if $\deg G=0$, $G$ is a constant that vanishes somewhere (such classes exist: if $N(Y)\subseteq\{0\}$, the first coordinate function vanishes on $N(Y)$ but not at $(1,0,\dots,0)$, contrary to (i)), so $G=0$. By~(i), $G(1,0,\dots,0)=0$,
i.e.\ $G$ vanishes at $[dx_1]=[1:0:\dots:0]$. As $G$ was arbitrary,
$[dx_1]\in\operatorname{Dir}(Y)$.
\end{proof}

\formalref{\leandecl{Stafford38/Paper/Section5/Proof.lean}{66}{Stafford38.Paper.Section5.theorem_5_7_i_paper};
\leandecl{Stafford38/Paper/Section5.lean}{51}{Stafford38.Paper.Section5.theorem_5_7_i};
\leandecl{Stafford38/Paper/Section5/Place.lean}{51}{Stafford38.Paper.Section5.lemma_5_9};
\leandecl{Stafford38/Paper/Section5/Place.lean}{155}{Stafford38.Paper.Section5.lemma_5_10};
\leandecl{Stafford38/Paper/Section5/ResidueBasis.lean}{33}{Stafford38.Paper.Section5.lemma_5_11};
\leandecl{Stafford38/Paper/Section5/KaehlerRankBound.lean}{139}{Stafford38.Paper.Section5.lemma_5_12_a_bound};
\leandecl{Stafford38/Paper/Section5/TangentCone.lean}{100}{Stafford38.Paper.Section5.lemma_5_12_a_cone};
\leandecl{Stafford38/Paper/Section5/TangentCone.lean}{142}{Stafford38.Paper.Section5.lemma_5_12_b};
\leandecl{Stafford38/Paper/Section5/TangentCone.lean}{171}{Stafford38.Paper.Section5.lemma_5_13}}

\subsection{Proof of Theorem \ref{thm:coisotropic-exclusion}}\label{sub:excl}
\begin{proof}
Suppose we have a \textit{nonempty} closed subset $W\subseteq T^*X$ that is fibre-conical and involutive. Suppose further, that $W\cap\pi^{-1}(H)=\varnothing$. 

We will give a proof by contradiction, with this rough logic: Suppose $W\neq\varnothing$ and let $Y = \pi(W)$. Since $W$ is by assumption fibre-conical, we start by showing that $Y\times \{0\} = \pi(W)\times \{0\} \subseteq W$ and $Y \cap H = \varnothing$.  Using involutivity, we then show that we can translate along $W$ in the direction of the fibres which means that all conormal covectors of every irreducible component $Y_0$ of $Y$ lie in $W$. This means $P$ vanishes on all of $W$, hence on $N(Y_0)$. This gives us the setup of Theorem \ref{thm:asymptotic-conormal-formal}, which implies $P(1,0,\dots,0) = 0,$ a contradiction to our assumption. \\

We define the (radical) vanishing ideal $J=I(W)\subseteq K[x,\xi]$ consisting of all polynomials that vanish at every point $(y,\xi)\in W$, and $B=J\cap K[x]$, consisting of the functions $f(x)$ with $f(y)=0$ for every $(y,\xi) \in W$. \\
Because $W$ is fibre-conical, the $\xi$-homogeneous components of any
polynomial in $J$ again lie in $J$: if $g=\sum_dg_d\in J$ and
$(y,\xi)\in W$, then $\sum_da^dg_d(y,\xi)=g(y,a\xi)=0$ for all $a\in K^\times$,
and a polynomial in $a$ with infinitely many roots is zero, so every
$g_d(y,\xi)=0$ and therefore $g_d \in J$. 
Let $g \in J$ be written in its $\xi$-homogeneous components form. In particular, $g_0 \in J$. Since $g_0$ has $\xi$-degree $0$, it is a polynomial in $x$ alone, so $g_0 \in B$. Let $Y = V_K(B)$. For $y \in Y$ we therefore have $g(y, 0) = g_0(y) = 0$, so every $g \in J$ vanishes at $(y,0)$, i.e. $(y,0)\in V_K(J)$. $V_K(J) = V_K(I(W)) = W$, since $W$ is closed \cite[Prop. I.1.2]{hartshorne},  hence $(y, 0) \in W$, or $Y \times \{0\} \subseteq W$. 
Furthermore, $Y$ is the projection of $W$ and is therefore per definition of $Y$ closed: $\pi(W)\subseteq Y$ since for ($y',\xi')\in W$ and $f \in B\subseteq J$, we have $f(y') = f(y',\xi') = 0$ so $y' \in V_K(B) = Y$, and $Y\subseteq \pi(W)$, since if $y \in Y$, then also $(y,0)\in W$, so $y = \pi(y,0) \in \pi(W)$.\\
Lastly, assume there exists a $y \in Y\cap H$, then $(y,0)\in W$ and $\pi(y,0) = y \in H$, which would put $(y,0)\in W\cap \pi^{-1}(H)$, which is wrong by our assumption. Therefore, $Y \cap H = \varnothing$.\\

For the next argument, we also refer to \cite[Ch. E.2]{HTT}. Let $f\in B$. Then $f\circ\pi\in J$.
The Hamiltonian derivation of $f$ applied to $g \in J$ is
\[ H_f g=\sum_i\frac{\partial f}{\partial x_i}
 \frac{\partial g}{\partial\xi_i} = \{f,g\}\]
since the second term of the Poisson bracket vanishes because $f$ does not depend on $\xi$. Since $f, g \in J$ and $\{J, J\} \subseteq J$ by the involutivity of $W$, it follows that $H_f g \in J$ and by induction also $H^m_fg \in J$. 
For a fixed $y \in Y$ and $g \in J$, we define the polynomial $\varphi(t)=g(y,t\,df_y)\in K[t]$. Its derivative is $\varphi'(t) = \sum_i\partial_{x_i}f(y)\frac{\partial g}{\partial \xi_i}{(y, tdf_y)} = H_fg(y, tdf_y)$, so by iteration, $\varphi^{(m)}(t) = H^m_fg(y, tdf_y)$ for $m\geq 0.$ 
Evaluating at $t = 0$, $\varphi^{(m)}(0) = H_f^mg(y,0).$ Since we have established that $(y,0)\in W$ and $H^m_fg \in J$, and all polynomials in $J$ vanish on $W$, we have $\varphi^{(m)}(0) = 0$ for all $m\geq 0$ -- so $\varphi(t)$ is a polynomial in $t$ all of whose derivatives vanish at $t = 0$.  If we write $\varphi(t) = \sum_m c_mt^m$, then $\varphi^{(m)}(0) = m!c_m$ and since $m!\neq 0$ because $K$ has characteristic zero, $c_m = 0$ for all $m$, and therefore, $\varphi(t) = 0$ everywhere. Remembering the definition of $\phi = g(y,t\,df_y) = 0$, so $(y,t\,df_y) \in V_K(J) = W$.
More generally, since any linear combination of elements in $B$ is again in $B$, we can apply this result to $f = \sum_a c_af_a \in B$ and evaluating at $t = 1$, it follows that 
\[(y, df_y) = (y, d(\sum_ac_af_a)_y) = (y, \sum_a c_ad(f_a)_y) \in W,\]
so the whole linear span $\{y\}\times \operatorname{span}\{df_y:f \in B\}$ is a subset of $W$.\\

Since $V_K(B)=Y=\pi(W)\neq\varnothing$, $B \neq K[x]$, so $B$ must have at least one minimal prime. The ideal $B$ is radical (by construction since it is a vanishing ideal) in the Noetherian ring $K[x]$, so $B$ has only finitely many minimal primes $\mathfrak q_1,\dots,\mathfrak q_r$ and $B=\bigcap_i\mathfrak q_i$ \cite[Cor. I.1.4 and Prop. I.1.5]{hartshorne}. We choose $\mathfrak q=\mathfrak q_1$ and set $Y_0=V_K(\mathfrak q)$. Since $\mathfrak q$ is prime, $1\notin\mathfrak q$, so $Y_0\neq\varnothing$ by Theorem~\ref{thm:nss}; moreover $Y_0$ is closed and irreducible since $\mathfrak q$ is prime \cite[Cor.~I.1.4]{hartshorne}, and since primeness implies $\mathfrak q$ is radical, $\mathfrak q = I(Y_0)$ by Theorem \ref{thm:nss}. Additionally, $Y_0\subseteq Y$, hence
$Y_0\cap H=\varnothing$, as we have already established $Y\cap H=\varnothing$.

Furthermore, if $r = 1$, we choose $s = 1$, then $sh \in \mathfrak q = B$, and if $r\geq 2$, choose for every $i\geq 2$ elements $s_i \in \mathfrak q_i \setminus \mathfrak q$. Then $s = s_2s_3\dots s_r$ does not lie in $\mathfrak q$ since no factor lies in $\mathfrak q,$ but does lie in every other $\mathfrak q_i$ since they are ideals and one of the factors lies in them. On the other hand, if $h \in \mathfrak q$, then $sh \in \mathfrak q$ since $h\in \mathfrak q$, and $sh \in \mathfrak q_i$ since $s \in \mathfrak q_i$ for all $i\geq 2$. Therefore, $sh \in \bigcap_i \mathfrak q_i = B,$ i.e. $s\mathfrak q \subseteq B.$ We now fix such a $h \in \mathfrak q.$ For a $y \in Y_0\cap\{s\neq0\}$, $(y, dh_y) \in W:$
Since $sh \in B,$ and $d(sh)_y = s(y)dh_y+h(y)ds_y = s(y)dh_y$ since $h(y) = 0$, so $(y, s(y)dh_y) \in W$ and since $s(y)\neq 0,$ also $(y, dh_y)\in W$. 
For $g\in J$, the polynomial $G(x)=g(x,dh_x)$ therefore vanishes on $Y_0\cap\{s\neq0\}$, and $sG$ vanishes on all of $Y_0$, so $sG \in \mathfrak q$ and therefore $G \in \mathfrak q$, which means that $G$ vanishes on all of $Y_0$. Since this is true for every $g \in J$, we have $g(y, dh_y) = G(y) = 0$ for all $g\in J$, so $(y,dh_y) \in V_K(J)= W.$  As before, we can argue that since $\mathfrak q$ is an ideal, it is closed under $K$-linear combinations, so for $h_1,\dots ,h_k \in \mathfrak q$ and $c_a \in K$, also $(y, \sum_a c_a (dh_a)_y) \in W$, or put differently, $\{y\} \times \operatorname{span}\{dh_y: h \in \mathfrak q\} \subseteq W$ for all $y \in Y_0$. At a smooth point $y$ of $Y_0$, by Definition \ref{def:conormal}, the conormal space is 
\[N_y^*Y_0 = \operatorname{span}\{df_y : f \in \mathfrak q\}\]
and we have already established that $\mathfrak q = I(Y_0)$. We therefore have that $\{y\} \times N_y^*Y_0 \subseteq  W$ for every smooth $y \in Y_0$, so $P$, which only depends on $\xi$, vanishes on $N(Y_0)$. Because we have established, that $Y_0$ is closed, irreducible, and $Y_0 \cap H = \varnothing $, $Y_0$ therefore fulfills all requirements of Theorem \ref{thm:asymptotic-conormal-formal}, so $P(1,0,\dots,0) = 0,$ contradicting our initial assumption. Therefore, $W$ must be empty.
\end{proof}

\formalref{%
\leandecl{Stafford38/Paper/Section5/Exclusion.lean}{447}{Stafford38.Paper.Section5.theorem_5_8};
\leandecl{Stafford38/Paper/Section5/SmoothLocus.lean}{645}{Stafford38.Paper.Section5.theorem_5_7_i_printed}}

\section{Proof of Theorem~\ref{thm:main}}\label{sec:proof}

\formalref{\leandecl{Stafford38/Paper/Section3/OrderSymbol.lean}{290}{Stafford38.Paper.Section3.proposition_3_6}}
\begin{proposition}\label{prop:char-avoids}
$\Char(Q)\cap V(x_1)=\varnothing,$
where $V(x_1)\subseteq \Spec k[x, \xi]$ is the set of primes containing $x_1$.
\end{proposition}
\begin{proof}
Write $M=\gr Q$. Since $M$ is finite over $k[x,\xi]$,
$\Supp_{k[x,\xi]}(M/x_1M)=\Supp_{k[x,\xi]}(M)\cap V(x_1)$ \cite[Tag~00L3]{Stacks}. By Corollary~\ref{cor:Uzero},
$M/x_1M=0$, so its support is empty, hence $\Supp_{k[x,\xi]}(M)\cap V(x_1)=\varnothing$.
Since $M$ is finitely generated, $\Char(Q)=V(\Ann_{k[x,\xi]}M)=\Supp_{k[x,\xi]}(M)$
\cite[Tag~00L2]{Stacks}. Therefore $\Char(Q)\cap V(x_1)=\varnothing$.
\end{proof}

\formalref{%
\leandecl{Stafford38/Paper/Section6/Avoidance.lean}{71}{Stafford38.Paper.Section6.proposition_6_1}}

\begin{lemma}\label{lem:base-change}
Let $K/k$ be a field extension and put
$\mathcal I_K=\tilde dA_n(K)+x_1^N\tilde dA_n(K)$, $Q_K=Q\otimes_kK=A_n(K)/\mathcal I_K$, with the
order filtration $Q_{K,i}=q\,F_iA_n(K)$. Then
\[
 \gr Q_K\cong\gr Q\otimes_kK,\qquad
 \Ann_{K[x,\xi]}\gr Q_K=\bigl(\Ann_{k[x,\xi]}\gr Q\bigr)\,K[x,\xi],
\]
and $Q=0\iff Q_K=0$. In particular, $\Char(Q_K)$ is the preimage of $\Char(Q)$
under $\Spec K[x,\xi]\to\Spec k[x,\xi]$.
\end{lemma}
\begin{proof}
By the PBW basis \cite[Prop 1.2.1]{coutinho}, $F_iA_n(K)=F_iA_n(k)\otimes_kK$ and $\mathcal I_K=\mathcal I\otimes_kK$, so
$Q_{K,i}=Q_i\otimes_kK$. As $K$ is flat over $k$,
$Q_{K,i}/Q_{K,i-1}=(Q_i/Q_{i-1})\otimes_kK$, and summing over $i$ gives
$\gr Q_K=\gr Q\otimes_kK$. The module $\gr Q$ is finitely generated (cyclic,
generated by the class of $q$) over the Noetherian ring $k[x,\xi]$, and
$k[x,\xi]\to K[x,\xi]$ is flat, so the annihilator of the finite module
commutes with this base change \cite[Tag~07T8]{Stacks}. Finally, $Q\otimes_kK=0$
implies $Q=0$, because $K$ is a nonzero free $k$-module.
\end{proof}

Let $d\in A_n(k)$ be nonzero of Bernstein degree $N$ and let
$\tilde d=c^{-1}\Phi_g(d)$ be its monic normal form. If $n=0$, then $A_0(k)=k$ and every nonzero $d$ is
a unit anyways, so let $n\ge1$. As in Section \ref{sec:euler} we write
\[ \mathcal I=\tilde dA_n(k)+x_1^N\tilde dA_n(k),\qquad Q=A_n(k)/\mathcal I .\]
It suffices to show $Q=0$: then $1\in \mathcal I$, i.e.\ there are
$R,S\in A_n(k)$ with $1=\tilde dR+x_1^N\tilde dS$.

Let $K=\bar k$, and let
\[W = V_K\bigl(\operatorname{Ann}_{K[x,\xi]}\operatorname{gr}Q_K\bigr)\subseteq K^n\times (K^n)^*,\]
be the set of closed points of $\operatorname{Char}(Q_K)$ (Theorem~\ref{thm:nss}),
with vanishing ideal $I(W)=\sqrt{\operatorname{Ann}\operatorname{gr}Q_K}$.
We check
the hypotheses of Theorem~\ref{thm:coisotropic-exclusion} for $W$ and the
polynomial $P=\sigma^{\mathrm{\ord}}_N(\tilde d)$ of Section~\ref{sec:monic}:

\begin{itemize}
\item $W$ is Zariski closed, since it is the zero set of an ideal of polynomials.
\item The module $\gr Q_K$ is graded by the order filtration, which assigns degree
$0$ to the $x_i$ and degree $1$ to the $\xi_i$; so $\Ann\gr Q_K$  is homogeneous in $\xi$. Hence if $(y,\xi)\in W$
and $g\in\Ann\gr Q_K$, then $g(y,t\xi)=\sum_dt^dg_d(y,\xi)=0$ for every
$t\in K^\times$, where $g_d$ denotes the $\xi$-homogeneous component of degree $d$ of $g$. So
$(y,t\xi)\in W$, i.e.\ $W$ is fibre-conical.
\item By Gabber's theorem (Theorem~\ref{ithm:gabber}) $I(W)$ is involutive. 
\item By eq.~\eqref{eq:P-annihilates}, $P\in
\Ann_{k[x,\xi]}\gr Q$, hence $P\in\Ann_{K[x,\xi]}\gr Q_K$ by
Lemma~\ref{lem:base-change}, so $P$ vanishes on $W$.
\item By eq.~\eqref{eq:P-e1}, $P(1,0,\dots,0)=1\neq 0$.
\item  By Lemma~\ref{lem:base-change}, $\operatorname{Char}(Q_K)$ is the preimage of
  $\operatorname{Char}(Q)$ under $\operatorname{Spec}K[x,\xi]\to\operatorname{Spec}k[x,\xi]$,
  and $V(x_1)\subseteq\operatorname{Spec}K[x,\xi]$ is the preimage of
  $V(x_1)\subseteq\operatorname{Spec}k[x,\xi]$. Since preimages commute with
  intersections, Proposition~\ref{prop:char-avoids} gives
  $\operatorname{Char}(Q_K)\cap V(x_1)=\emptyset$. By Theorem~\ref{thm:nss}(a) and Definition~\ref{def:spec}, a closed point $\mathfrak m_{(y,\xi)}$ lies in $V(x_1)$ if and only if $y_1=0$, i.e.\ the closed points of $V(x_1)$ form $\pi^{-1}(H)=H\times(K^n)^*$. Hence the closed points of $\operatorname{Char}(Q_K)\cap V(x_1)$ are $W\cap\pi^{-1}(H)$, and so $W\cap\pi^{-1}(H)=\emptyset$.
\end{itemize}
By Theorem~\ref{thm:coisotropic-exclusion}, $W=\varnothing$. As $K$ is
algebraically closed, by Theorem \ref{thm:nss}  $1 \in \Ann_{K[x,\xi]}\gr Q_K$,
so $\gr Q_K=0$. Since the filtration of $Q_K$ is exhaustive with $Q_{K,-1}=0$,
induction on $m$ implies $Q_{K,m}=0$ for all $m$, so $Q_K=0$, and
$Q=0$ by Lemma~\ref{lem:base-change}. Therefore, $A_n(k) = \tilde dA_n(k) + x_1^N\tilde dA_n(k)$, or equivalently $1=\tilde dR+x_1^N\tilde dS$
with $R,S\in A_n(k)$.

We apply the inverse of the automorphism $\Phi_g$ of
Remark~\ref{rem:transport}:
\[
 1=d\bigl(c^{-1}\Phi_g^{-1}(R)\bigr)+F\,d\bigl(c^{-1}\Phi_g^{-1}(S)\bigr),
 \qquad F=\Phi_g^{-1}(x_1)^N,
\]
which proves Theorem~\ref{thm:main}. \qed

\formalref{\leandecl{Stafford38/Paper/Main.lean}{24}{Stafford38.Paper.theorem_1_1},
\leandecl{Stafford38/Paper/Main.lean}{29}{Stafford38.Paper.theorem_1_1_fixedSource},
\leandecl{Stafford38/Paper/Section6/Proof.lean}{357}{Stafford38.Paper.Section6.theorem_1_1_paper},
\leandecl{Stafford38/Paper/Section6/Proof.lean}{318}{Stafford38.Paper.Section6.theorem_6_Q_eq_zero}}

\formalref{%
\leandecl{Stafford38/Paper/Section6/BaseChange.lean}{460}{Stafford38.Paper.Section6.lemma_6_2}}

Lastly, as Stafford already showed, from $1 = dR+FdS$ follows cyclicity of the torsion right module:
\begin{corollary}[\cite{Stafford1978}]\label{cor:cyclic}
Every finitely generated torsion right $A_n(k)$-module $T$ is cyclic, i.e. can be generated by one element. Moreover, it is stably cyclic, meaning that given any $l\geq 1$ and $a_1,...a_{l+1} \in T$ that generate $T$, there exist $f_i \in A_n$ such that $(a_i + a_{l+1}f_i)$ $i = 1,...l$ also generate $T$ \cite[p. 435]{Stafford1978}. 
\end{corollary}

\begin{proof}
    Suppose for every $d\neq 0 \in A_n$, there exists $F, R, S \in A_n$, such that $1 = dR + FdS$. Let us denote a torsion module by $T$. Stafford's theorem (the right module analogue to Theorem 3.7 \mbox{\cite[p.~438]{Stafford1978}}) established, that every finitely generated (right) torsion module can be generated by 2 generators, so let us denote the generators of $T$ with $t_1$ and $t_2$. We can define a homomorphism
    \[\phi:A_n \oplus A_n\rightarrow T \quad \phi(r,s)\mapsto t_1r+t_2s.\]
    This map is surjective, since $t_1, t_2$ generate all of $T$. By the first isomorphism theorem (for $\varphi:M\rightarrow N$, $im(\varphi) \cong M/ker(\varphi)$), it follows
    \[T \cong (A_n\oplus A_n)/M\]
    with $M = ker(\phi)$. \\
    $T$ is torsion and therefore $\forall t\in T: \exists m \in A_n\backslash\{0\}: tm = 0.$ Specifically, let us denote with $d_1$ and $d_2$ the annihilators of $t_1$ and $t_2$, i.e.
    \[t_1d_1 = t_2d_2 = 0.\]
    $A_n$ is a right Ore domain \cite[Prop. 1.8.4]{Bjork}, which means it has the property that for any nonzero $d_1, d_2$ there exists elements $r_1, r_2 \in A_n$ such that $d_1r_1 = d_2r_2 \neq 0.$
    Defining therefore $d := d_1r_1 = d_2r_2$ yields us a common annihilator of $t_1$ and $t_2$.\\
    $(d, 0)$ and $(0,d)$ are elements in $M = ker(\phi)$, since $\phi(d, 0) = t_1d = 0$ and $\phi(0, d) = t_2d = 0.$\\
    Because $(0,d), (d, 0) \in M$, also
    $(d, 0)S + (0, d)R = (dS, dR) = (dS, 1-FdS) =: (a, 1-Fa) \in M$, i.e. $(a, 1-Fa)$ is an element of $M$. \\
    On the other hand, $A_n\oplus A_n$ can be generated by $(a, 1-Fa)$ and $(1, -F)$. We can check this by finding $u,v$ such that for all $(x,y)\in A_n\oplus A_n$ we can write
    \[(x,y) = (a, 1-Fa)u + (1, -F)v.\]
    Suitable choices are $u = y+Fx$ and $v = x-a(y+Fx)$.\\
    Since  $A_n\oplus A_n$ can be generated by $(a, 1-Fa)$ and $(1, -F)$ and $(a, 1-Fa)$ is an element of $M$, it follows that $T \cong (A_n\oplus A_n)/M$ can be generated by $\phi(1, -F)$, i.e. $T$ is cyclic. \\
    For stable cyclicity we see that since $t_1, t_2$ were arbitrary generators, for $l = 1$ we can set $a_1 = t_1, a_2 = t_2$ and choose $f_1 = -F$, then $\phi(1, -F) = t_1 - t_2F$ generates $T$. For $l\geq 2$ we simply note that by the right module version of \cite[Theorem~3.7]{Stafford1978}, $T$ is already proven to be stably two-generated, therefore, T is stably cyclic.
    
\end{proof}
\formalref{\leandecl{Stafford38/Paper/Main.lean}{67}{Stafford38.Paper.corollary_6_3},
\leandecl{Stafford38/Paper/Main.lean}{81}{Stafford38.Paper.corollary_6_3_stablyCyclic},
\leandecl{Stafford38/Paper/Section6/Cyclicity.lean}{219}{Stafford38.Paper.Section6.corollary_6_3_of};
Mathlib-only statements \leandecl{Stafford38/ChallengeDefinitions.lean}{302}{Stafford38Challenge.PaperCorollary_6_3},
\leandecl{Stafford38/PaperClaimDefinitions.lean}{33}{Stafford38Challenge.PaperCorollary_6_3_stablyCyclic}}

\section{Generalizations and open questions}\label{sec:generalizations}
The formalization in Lean discussed in Section~\ref{sec:formal} also proves the following consequences, which have no printed proof: the left-handed identity $1=Rd+SdF$, via the
formal adjoint, which implies Corollary \ref{cor:cyclic} also for left torsion modules; the identity for finite-order differential operators on localizations of $k[x_1,\dots,x_n]$, including principal opens,
Laurent-polynomial coefficients and the fraction field $k(x_1,\dots,x_n)$; and, independently
of the geometric part of Section \ref{sec:conormal-review}, the identity with $F=x^r$ for
$d=p^r-\sum_j a_jx^{m_j}$, $r\ge1$, in any $k$-algebra containing elements
$x,p$ with $px=xp+1$, provided each $a_j$ commutes with $x$ and $p$.
\formalref{\leandecl{Stafford38/Paper/Section7.lean}{38}{Stafford38.Paper.Section7.leftHandedIdentity};
\leandecl{Stafford38/Paper/Section7.lean}{72}{Stafford38.Paper.Section7.corollary_6_3_left};
\leandecl{Stafford38/Paper/Section7.lean}{225}{Stafford38.Paper.Section7.localizedDifferentialIdentity};
\leandecl{Stafford38/Paper/Section7.lean}{246}{Stafford38.Paper.Section7.evolutionaryIdentity}}

Two further results, also proved in Lean, give the identity $1=dR+FdS$ in larger classes of rings, though without the
choice of $F$ in Theorem~\ref{thm:main}: Let $k$ be a field of
characteristic zero.
\begin{enumerate}[label=(\alph*),nosep,leftmargin=1.6em]
\item Let $A$ be a smooth finitely generated $k$-algebra that is a domain,
and let $D_k(A)$ be the ring of $k$-linear differential operators of finite
order on $A$. For every nonzero $d\in D_k(A)$ there are $F,R,S\in D_k(A)$
with $1=dR+FdS$~\cite{GlobalStaffordPalomar2026,GlobalStaffordZenodo106}.
For $A=k[x_1,\dots,x_n]$ this is the identity in $A_n(k)$.
\formalref{\leanglobal{Solution.lean}{18}{GlobalStaffordChallenge.universalStatement};
statement \leanglobal{Challenge.lean}{70}{GlobalStaffordChallenge.UniversalStatement}}
\item Let $U$ be an infinite-dimensional simple $k$-algebra with an
exhaustive filtration $k=F_0U\subseteq F_1U\subseteq\cdots$ by
finite-dimensional subspaces, $F_iU\cdot F_jU\subseteq F_{i+j}U$, whose
associated graded ring is a commutative domain finitely generated over $k$.
For all nonzero $d,a\in U$ and every $B\in U$ there are $H,R,S\in U$ with
$1=dR+(B+Ha)dS$~\cite{SameDivisorFiltered}. The Bernstein filtration makes
$A_n(k)$ such an algebra, and $B=0$, $a=1$ gives $1=dR+HdS$.
\formalref{\leanfiltered{Solution.lean}{46}{SymplecticStaffordChallenge.universalStatement};
statement \leanfiltered{Challenge.lean}{54}{SymplecticStaffordChallenge.AffineSameDivisorStatement}}
\end{enumerate}
The proof of Theorem~\ref{thm:main} shows $Q=0$ without producing $R$ and
$S$. Over a computable field, the fixed choice of $F$ nevertheless gives a terminating algorithm:
enumerate the countably many pairs $(g,D)$ of linear symplectic changes of generators and degree
bounds, and solve the linear system $1=dR+FdS$ with
$F=\Phi_g^{-1}(x_1)^{\deg_Bd}$ for $R,S$ of Bernstein degree at most $D$.
Theorem~\ref{thm:main} guarantees that some pair succeeds. We cannot provide bounds for
$\deg_BR$ and $\deg_BS$ in terms of $\deg_Bd$, so the question for an efficient
algorithm remains open. In (a) and (b), $F$ is not fixed, so this search does
not apply.

\backmatter

\section*{Declarations}

\bmhead{Acknowledgements}
J. Moser was funded by the Austrian Research Promotion Agency (FFG) within
the programme AI for Green, project VENTUS (910263). The
authors would like to acknowledge the use of HPC resources provided by the ZID
of Graz University of Technology. We gratefully acknowledge support from NAWI
Graz.

\bmhead{Competing interests}
The authors have no competing interests to declare that are relevant to the
content of this article.

\bmhead{Code availability}
The Lean sources of the formal proof and the supplementary material are
archived~\cite{Stafford38Zenodo140,Stafford38Supplementary033}.

\bibliography{references}

\begin{appendices}

\section{AI-assisted research and formalization}\label{app:ai-workflow}

AI agents assisted with proof discovery, counterexample searches,
review, formalization and writing. The work began with computations of torsion generators in
Macaulay2 and Maple. Direct attempts to construct universal cofactors
produced special cases and reductions, but no general construction. The
decisive step was to study the quotient $Q=A_n/(dA_n+x_1^N dA_n)$ and to show
that its characteristic support is empty.

The first complete argument was written by AI agents in the style of a paper
proof. Agents then formalized it in Lean, and the formalization changed parts
of the argument where the original steps could not be checked as written.
After the Lean proof was complete and published, the paper proof was rewritten
from the checked Lean proof, so that every statement corresponds to a Lean
declaration, and the Lean proof was then again aligned with the printed proof
step by step. The authors reviewed this correspondence and substantially revised the text for
clarity and readability together with the agents, until neither the human review nor the AI
checks raised further objections.

For finding the first complete argument, we proceeded as follows:
One controlling agent kept the authoritative state of the proof.
Parallel agents worked on separate approaches and returned evidence together
with its scope, and independent reviewer agents checked each result before it was
accepted. The first author steered the work with short instructions, and the
agents wrote the detailed task prompts for each other. Advice from
consultations counted only after an independent check.

The working rules developed together with the proof into a reusable
agent skill for open mathematical problems. Early versions emphasized
falsification and keeping conditional results apart from established ones.
Once a candidate proof existed, the emphasis moved to repairing the first
invalid step of a promising argument, and
to auditing proofs.

Research and review used OpenAI GPT models from GPT-5.5 to GPT-6.1,
Anthropic Claude models from Haiku~4.5 to Opus~5, and locally hosted Qwen
and DeepSeek models. GPT and Claude models, most recently Claude Opus~5.5,
helped draft and revise the manuscript and the supplementary material. The
authors are responsible for the final text.

The formal proof uses Lean~4, Mathlib and AlgebraicAnalysis
v0.3.3~\cite{Lean4,Mathlib2020,AlgebraicAnalysis033}. Conditional results
were kept apart from proved ones until every hypothesis was discharged.
Comparator checked that the proved theorems have exactly the statements of
separate challenge files that use only Mathlib definitions, and both the Lean
kernel and the independent checker NanoDa accepted the proofs
\cite{Comparator,NanoDa,PalomarDocs}. Theorem~\ref{thm:main} and
generalization~(a) of Section~\ref{sec:generalizations} are registered in the
Palomar registry~\cite{Stafford38Palomar2026,GlobalStaffordPalomar2026}, and
their Lean sources are archived~\cite{Stafford38Zenodo140,GlobalStaffordZenodo106}.

\end{appendices}

\end{document}